\documentclass[amstex]{article}
\usepackage{graphicx}
\usepackage{dcolumn}
\usepackage{amsmath}
\usepackage{amssymb}
\usepackage{amsfonts}
\usepackage{amscd}
\begin{document}
\newtheorem{thm}{Theorem}[section]
\newtheorem{lem}{Lemma}[section]
\newtheorem{prop}{Proposition}[section]
\newtheorem{cor}{Corollary}[section]
\newtheorem{assum}{Assumption}[section]
\newtheorem{rem}{Remark}[section]
\newtheorem{defn}{Definition}[section]
\newcommand{\lv}{\left \vert}
\newcommand{\rv}{\right \vert}
\newcommand{\lV}{\left \Vert}
\newcommand{\rV}{\right \Vert}
\newcommand{\la}{\left \langle}
\newcommand{\ra}{\right \rangle}
\newcommand{\ltm}{\left \{}
\newcommand{\rtm}{\right \}}
\newcommand{\lcm}{\left [}
\newcommand{\rcm}{\right ]}
\newcommand{\ket}[1]{\lv #1 \ra}
\newcommand{\bra}[1]{\la #1 \rv}
\newcommand{\lmk}{\left (}
\newcommand{\rmk}{\right )}
\newcommand{\al}{{\mathcal A}}
\newcommand{\md}{M_d({\mathbb C})}
\newcommand{\ali}[1]{{\mathfrak A}_{[ #1 ,\infty)}}
\newcommand{\alm}[1]{{\mathfrak A}_{(-\infty, #1 ]}}
\newcommand{\nn}[1]{\lV #1 \rV}
\newcommand{\br}{{\mathbb R}}
\newcommand{\dm}{{\rm dom}\mu}
\title{Approximating macroscopic observables in quantum spin systems
with commuting matrices}
\author{Yoshiko Ogata\thanks
{Graduate School of Mathematics, University of Tokyo, Japan}}
\maketitle{}
\begin{abstract}
Macroscopic observables in a quantum spin system
are given by sequences of spatial means of local elements
$\frac{1}{2n+1}\sum_{j=-n}^n\gamma_j(A_{i}),
\; n\in{\mathbb N},\; i=1,\cdots,m$
in a UHF algebra.
One of their properties is that they commute asymptotically, as $n$
goes to infinity.
It is not true that any given set of asymptotically commuting
matrices can be approximated by
commuting ones in the norm topology.
In this paper, we show that for macroscopic observables,
this is true.
\end{abstract}
%%%%%%%%%%%%%%%%%%%%%%%%%%%%%%%%%%
\section{Introduction}
The infinite quantum spin chain with one site algebra 
${M_d}({\mathbb C})$ is given by the UHF $C^*$-algebra
\[
{\mathcal A}_{{\mathbb Z}}:=\overline{
\bigotimes_{{\mathbb Z}}{M_d}({\mathbb C})}^{C^*},
\]
which is the $C^*$- inductive limit of the local algebras
\[
\left\{
{\mathcal A}_{\Lambda}:=\bigotimes_{{\Lambda}}{M_d}({\mathbb C})\vert
\quad \Lambda\subset{\mathbb Z},\quad
| \Lambda|<\infty
\right\}
.
\]
Here, $| \Lambda|$ denotes the number of points in $\Lambda$.
We denote ${\mathcal A}_{[-n,n]}$ by ${\mathcal A}_{n}$.
Let $\gamma_j,\; j\in {\mathbb Z}$ be the $j$-lattice translation.
We say that a state $\omega$ is translation invariant if $\omega\circ\gamma_j=\omega$
for all $j\in{\mathbb Z}$.
%A state $\omega$ is called translation-invariant if 
%$\omega\circ\gamma_j=\omega$ for all $j\in{\mathbb Z}$.
We denote the set of all translation invariant states
by $S_\gamma$.
For each finite subset $\Lambda$ of $\mathbb Z$, we set
${\rm Tr}_{\Lambda}$ to be the non-normalized trace over $\otimes_{\Lambda}M_d({\mathbb C})$. Furthermore, we denote ${\rm Tr}_{[-n,n]}$ by ${\rm Tr}_n$.
\begin{thm}\label{main}
Let $A_1,\cdots,A_m$ be self-adjoint elements 
in ${M_d}({\mathbb C})$, and define
\[
H_{i,n}:=\frac{1}{2n+1}\sum_{j=-n}^n\gamma_j(A_{i})
\in{\mathcal A}_n,
\] for each $i=1,\cdots,m,\; n\in{\mathbb N}$.
Then there exist sequences of self-adjoint elements
$Y_{i,n},\;i=1,\cdots, m,\;n\in{\mathbb N}$ such that
\begin{align*}
Y_{i,n}\in{\mathcal A}_n,\quad
\lim_{n\to\infty}\lV
H_{i,n}-Y_{i,n}
\rV=0,\quad \forall i=1,\cdots,m,\end{align*}
and\begin{align*}
\left[Y_{i,n},Y_{j,n}\right]=0,\quad\forall i,j=1,\cdots,m.
\end{align*}
\end{thm}
Following the standard procedure in statistical mechanics,
we can extend this theorem to general
translation invariant interactions.
(See Appendix \ref{gen})\\\\
From the Theorem of Lin \cite{lin1}, 
we already know that any pair of sequences of matrices whose commutator
vanishes asymptotically can be approximated by
commuting ones in general. Therefore, $m=2$ case is already known.
On the other hand, for $m\ge 3$, it is known 
that such a statement is not true in general \cite{dav}.
However, for macroscopic observables,
because of their nice thermodynamic structure,
we can show that the Theorem is true.\\\\
This article is organized as follows:
In section \ref{ef}, we introduce the entropy function
$\mu$ associated with $A_1,\cdots,A_m$.
This function gives an estimate for rank
of projections which are concentrating at value
$x=(x_1,\cdots,x_m)$ with respect to $A_1,\cdots,A_m$
(Lemma \ref{lem9}, Lemma \ref{lem10}).
In section \ref{tower},
we consider a $C^*$-algebra $\prod_{k}M_{n_k}({\mathbb C})/
\oplus_k M_{n_k}({\mathbb C})$.
For each $s\in{\mathbb R}$,
we define an ideal $I_s$ and investigate its property.
%This ideal $I_s$ is related to the rank of projections
%above.
We construct a tower of ideals, using the information of $\mu$
and $I_s$.
By virtue of this ideal tower,
in section \ref{proof},
we prove Theorem \ref{main}, using 
the technique developed in \cite{eglp},\cite{lin1}.
\section{The entropy function $\mu(x)$}\label{ef}
In the rest of this paper we
 fix self-adjoint elements $A_1,\cdots,A_m\in \md$.
For a translation invariant state 
$\omega\in S_\gamma({\al})$,
the mean entropy 
\begin{align*}
s(\omega):=\lim_{n\to\infty}\frac{1}{2n+1}
S_{[-n,n]}(\omega)
\end{align*}is known to exists (See\cite{BR96}).
Here, $S_{[-n,n]}(\omega)$ is the 
von Neumann entropy 
of $\omega\vert_{\al_{n}}$, the restriction
of $\omega$ to $\al_{n}$.
The function $S_{\gamma}({\al})\ni\omega\mapsto s(\omega)\in {\br}$
is affine and upper semi-continuous,
when $S_{\gamma}({\al})$ is equipped with
the weak$*$-topology.
Furthermore, it takes values in 
$[0,\log d]$.
\begin{defn}
The entropy function $\mu:{\mathbb R}^m
\to[-\infty,+\infty)$
 associated with
$A_1,\cdots,A_m$ is defined by
\begin{align}
\mu(x):=\left\{
\begin{gathered}
\sup\left\{
s(\omega):\omega\in S_\gamma({\al}),\;
\omega(A_i)=x_i,\; 1\le i\le m
\right\},\\
\quad{\rm if}\;
\{\omega\in S_\gamma({\al}):\;
\omega(A_i)=x_i,\; 1\le i\le m\}\neq\phi\\
-\infty,\\\quad {\rm if}\;
\{\omega\in S_\gamma({\al}):\;
\omega(A_i)=x_i,\; 1\le i\le m\}=\phi
\end{gathered}\right..
\end{align}
\end{defn}
By the definition, $\mu$ is concave, 
upper semi continuous and $\mu(x)\le {\rm log} d$ for all 
$x\in {\mathbb R}^m$.
We set the domain of $\mu$ by 
${\rm dom}\mu:=\{x\in {\mathbb R}^m \; :\; \mu(x)>-\infty\}$.
As the set $\{(\omega(A_1),\cdots,\omega(A_m))\; :\;
\omega\in S_{\gamma}({\al})\}$ is in 
$\prod_{i=1}^m[-\lV A_i\rV,\lV A_i\rV]$,
${\rm dom}\mu$ is bounded.
\\
We denote the level sets of $\mu$ by
\begin{align}
X_s:=\{x\in\br^m\;:\; \mu(x)\ge s\},\quad
s\in{\mathbb R}.
\end{align}
From the upper semi-continuity and 
the concavity of $\mu$,
$X_s$ is compact and convex.
Note that if $x\in dom\mu$
then $\mu(x)\in[0,\log d]$.
From this, we have
${\rm dom} \mu= X_0$.
%Clearly, ${\rm dom} \mu=X_0$ is a bounded set in $\br^m$,
%and closed because of the upper semi continuity
%of $\mu$. 
Therefore, ${\rm dom} \mu$ is compact.\\
The entropy function $\mu$ is the Legendre
transform of the {\it free energy} function $p$:
\begin{lem}\label{pmu}
Let $p:{\mathbb R}^m\to\br$ be a function defined by
\begin{align*}
p(\alpha):=\log {\rm Tr}_{0} e^{\sum_{i=1}^m\alpha_i A_i},
\quad \alpha\in{\br}^m.
\end{align*}
Then we have
\begin{align}\label{pmue}
 p(\alpha)=
\sup\left\{(\alpha,x)+
\mu(x)\;:x\in{\br}^m
\right\},\quad \forall \alpha\in{\br}^m,
\end{align} 
and
\begin{align}\label{mupe}
\mu(x)=-\sup\left\{
(\alpha,x)-p(\alpha)\;:\;\alpha\in\br^m
\right\}.
\end{align}
Here, $(\alpha,x)$ is the inner product of $\br^m$ : 
$(\alpha,x):=\sum_{i=1}^m\alpha_i x_i$.
\end{lem}
{\it Proof} See Appendix.
$\square$\\
Later, we will need contour lines of $\mu$
which are $\varepsilon$-dense in ${\rm dom}\mu$
in the following sense.
For $\varepsilon>0$ and a set $A\subset\br^m$,
we denote the $\varepsilon$-neighborhood of $A$ by $B_{\varepsilon}(A)$.
\begin{lem}\label{lem21}
For any $\varepsilon>0$,
there exists a finite sequence of
real numbers $s_0>s_1>\cdots>s_n$ with $s_n<0$,
such that
\begin{align*}
s_0=\sup\left\{
\mu(x)\;:\; x\in{\br}^m
\right\}<\infty,\\
X_{s_k}\subset B_\varepsilon(X_{s_{k-1}}),\quad k=1,\cdots,n,\\
{\rm dom}\mu=X_0= X_{s_n}.
\end{align*}
Furthermore, $X_{s_0}$ consists of one point
\[
x_0:=(\frac{{\rm Tr}_{0} A_1}{d},\cdots \frac{{\rm Tr}_{0} A_m}{d})\in\br^m.
\]
\end{lem}
{\it Proof}
Recall that $\mu$ takes values in $[-\infty,\log d]$,
and ${\rm dom}\mu$ is non-empty.
Therefore, $s_0:=\sup\left\{
\mu(x)\;:\; x\in{\br}^m
\right\}=\sup\left\{
\mu(x)\;:\; x\in
{{\rm dom} \mu}=X_0
\right\}$ is finite.
As ${{\rm dom} \mu}=X_0$ is compact
and $\mu$ is upper semi-continuous, 
there exists $x_0\in{\br}^m$ such that
$s_0:=\mu(x_0)
$.\\
To see that $X_{s_0}$ consists of one point $x_0$,
note from Lemma \ref{pmu} that 
$\mu(x)=s_0$ is equivalent to
$-p(0)\ge (\alpha,x)-p(\alpha)$ for all
$\alpha
\in{\br}^m$.
As $p$ is a differentiable function, 
this implies $x_i=\frac{\partial}{\partial\alpha_i}p(0)=\frac{{\rm Tr}_{0} A_i}{d}$.
Therefore $X_{s_0}$ consists of one point
$x_0:=(\frac{{\rm Tr}_{0} A_1}{d},\cdots \frac{{\rm Tr}_{0} A_m}{d})$.
\\
Fix $\varepsilon>0$.
We claim that there exists $s_1<s_0$
such that $X_{s_1}\subset B_\varepsilon(X_{s_0})$.:
Assume the claim were false.
Then for any increasing sequence of real numbers
$\{s_k\}_{k},\; s_k<s_0,\;s_k\uparrow s_0$, 
the sets $X_{s_k}\cap  B_\varepsilon(X_{s_0})^c$ are not empty.
Choose $x_k\in X_{s_k}\cap  B_\varepsilon(X_{s_0})^c$
for each $k$. By the compactness of $\dm$,
$\{x_k\}$ has a convergent subsequence $\{x_k'\}$,
$x_k'\to x\in{\br}^m$.
By the upper semi-continuity of $\mu$,
we have $x\in X_{s_0}$.
However, this means $x_k'$s are in $B_\varepsilon(X_{s_0})$
eventually, which is a contradiction.
Accordingly, we obtain the claim.
\\
Let $C$ be a finite positive number such that
$\max\{{\rm diam}( {\rm dom}\mu) ,\varepsilon\}<C$.
Our second claim is that
for any $s<s_0$, 
\begin{align}\label{ccc}
X_{s'}\subset B_{\varepsilon}(X_s),\quad
s-\frac{\varepsilon(s_0-s)}{C-
\varepsilon}\le
\forall s'<s
\end{align}
holds.
Let $s'$ be a real number such that $s-\frac{\varepsilon(s_0-s)}{C-
\varepsilon}\le s'<s<s_0$.
From the concavity of 
$\mu$, we have $(1-t)x+tx_0\in X_{(1-t)s'+ts_0}$
for any $x\in X_{s'}$ and 
$0\le t\le 1$.
As
\[
\lV
x-((1-t)x+tx_0)
\rV\le t\lV
x-x_0
\rV<tC,
\]
this implies
\[
X_{s'}\subset B_{tC}(X_{(1-t)s'+ts_0}),
\]
for any $t\in[0,1]$.
Note that $\frac{s-s'}{s_0-s'}C\le \varepsilon$,
if $s-\frac{\varepsilon(s_0-s)}{C-
\varepsilon}\le s'<s$.
Therefore setting $t:=\frac{s-s'}{s_0-s'}\in[0,1]$,
we obtain
$X_{s'}\subset B_{\frac{s-s'}{s_0-s'}C}(X_s)
\subset B_{\varepsilon}(X_s)$.\\
Now we define $s_k:=s_0-
\lmk\frac{C}{C-\varepsilon}\rmk^{k-1}\cdot(s_0-s_1),\; k\ge 1$.
These $s_k$ satisfy
$s_k=s_{k-1}-\frac{\varepsilon}{C-\varepsilon}\cdot
(s_0-s_{k-1})\le s_{k}<s_{k-1},\;k\ge 2$.
Therefore, by the above claim, we get 
$X_{s_k}\subset B_{\varepsilon}(X_{s_{k-1}})$.
As $\frac{C}{C-\varepsilon}>1$, 
there exists $n\in{\mathbb N}$ such that
$s_n<0$.
For this $n$, we have $\dm=X_0=X_{s_n}$.$\square$\\\quad\\
One important property of the entropy function $\mu$ is
that it gives an asymptotic estimate for rank of
projections.
We first give the upper bound:
\begin{lem}\label{lem9}
Let $C$ be a nonempty compact convex subset of $\br^m$, 
$\{n_k\}_{k=1}^{\infty}$
 a subsequence of $\mathbb N$, and $\{p_k\}_{k=1}^\infty$ 
a sequence of 
projections in $\al$ such that 
$p_{k}\in \al_{n_k},\; k\in{\mathbb N}$.
Suppose that for any $\varepsilon>0$, we have
\[
\lmk \frac{ {\rm Tr}_{n_k}p_{k}H_{1,n_k}}{ {\rm Tr}_{n_k}p_k},\cdots,
\frac{ {\rm Tr}_{n_k}p_kH_{m,n_k}}{{\rm Tr}_{n_k}p_k}
\rmk\in C_{\varepsilon},
\] 
eventually.
%, where $C_{\varepsilon}$ is an $\varepsilon$-neighborhood of $C$.
Then we have
\begin{align}\label{ub9}
\limsup_k\frac 1{2n_k+1} \log {\rm Tr}_{n_k} p_{k}\le
\sup\{\mu(x):x\in C\}.
\end{align}
\end{lem}
{\it Proof} 
First we claim 
\begin{align}\label{ubp}
\limsup_{k}\frac 1{2n_k+1} \log {\rm Tr}_{n_k} p_k
\le\inf\left\{
\sup\left\{
p(\alpha)-(\alpha,x)
:\quad x\in C
\right\}:\alpha\in\br^m \right\}.
\end{align}
To prove this, we use an argument in \cite{ee}.
Fix $\alpha\in\br^m$.
By the positivity of the relative entropy, we have
\[
0\le S(\frac{ p_k}{ {\rm Tr}_{n_k}p_k}, 
\frac{e^{(2n_k+1)\sum_{i=1}^m
\alpha_i\cdot H_{i,n_k}}}{{\rm Tr}_{n_k} e^{(2n_k+1)\sum_{i=1}^m
\alpha_i\cdot H_{i,n_k}}}).\]
From this, we obtain 
\[
\frac 1{2n_k+1} \log {\rm Tr}_{n_k} p_k\le
p(\alpha)-
\sum_{i=1}^m\alpha_i \frac{{\rm Tr}_{n_k} p_{k} H_{i,n_k}}{{\rm Tr}_{n_k} p_k}.
\]
By the assumption, for any $\varepsilon>0$, we have
\[
\sum_{i=1}^m\alpha_i \frac{{\rm Tr}_{n_k} p_k H_{i,n_k}}{{\rm Tr}_{n_k} p_{k}}
\ge \inf\{(\alpha,x) :x\in C_{\varepsilon}\},\] eventually.
Therefore, we get
\begin{align*}
\limsup_{k}\frac 1{2n_k+1} \log {\rm Tr}_{n_k} p_k
\le
p(\alpha)-
\inf\left\{
(\alpha,x)
:\quad x\in C_{\varepsilon}
\right\},
\end{align*}
for all $\varepsilon>0$.
Taking $\varepsilon \to 0$ limit, we obtain 
\begin{align*}
\limsup_{k}\frac 1{2n_k+1} \log {\rm Tr}_{n_k} p_k
\le
p(\alpha)-
\inf\left\{
(\alpha,x)
:\quad x\in C
\right\}
=\sup\left\{p(\alpha)-
(\alpha,x)
:\quad x\in C
\right\},
\end{align*}
for all $\alpha\in\br^m$.
From this we have
\begin{align}\label{fs}
\limsup_{k}\frac 1{2n_k+1} \log {\rm Tr}_{n_k} p_k
\le
\inf\left\{
\sup\left\{p(\alpha)-
(\alpha,x)
:\quad x\in C
\right\}
\; :\; \alpha\in\br^m
\right\}.
\end{align}
The last term in (\ref{fs}) can be written as 
\begin{align}
\inf\left\{
\sup\left\{p(\alpha)-
(\alpha,x)
:\quad x\in C
\right\}
\; :\; \alpha\in\br^m
\right\}
=-\sup\left\{\inf\left \{
(\alpha,x)-p(\alpha):\quad x\in C\right\}\alpha\in{\br^m}\right\}.
\end{align}
From Sion's lemma \cite{sion}, we have
\begin{align*}
\sup\left\{\inf\left \{
(\alpha,x)-p(\alpha):\quad x\in C\right\}\alpha\in{\br^m}\right\}
=\inf\left\{\sup\left \{
(\alpha,x)-p(\alpha):\quad \alpha\in{\br^m}\right\}x\in C\right\}.
\end{align*}
By the equality (\ref{mupe}), the last term is equal to
 $\inf\left\{-\mu(x)\; :\;x\in C\right\}=-\sup
\{\mu(x)\; : x\in C\}$.
Combining this and (\ref{fs}), we obtain (\ref{ub9}).
$\square$\\\quad\\
In order to prove the lower bound, we use the following Theorem
from \cite{qsm} : 
\begin{thm}[\cite{qsm}]\label{Bj}
Let $\omega$ be an ergodic state over $\al$
and define
\begin{align*}
\beta_{\varepsilon,n}(\omega)
:=\min \left\{
\log {\rm Tr}_{n} q\; : \; q\in {\rm Proj}({\al}_n),\; \omega(q)\ge 1-\varepsilon
\right\}, 
\end{align*}
for each $0<\varepsilon<1$.
Then we have
\begin{align*}
\lim_{n\to\infty}\frac{1}{2n+1}
\beta_{\varepsilon,n}(\omega)=s(\omega).
\end{align*}
\end{thm}
The lower bound is given as follows:
\begin{lem}\label{lem10}
Let $U_1,\cdots,U_m$ be open subsets of $\br$,
$\{n_k\}$  a subsequence of $\mathbb N$,
and $\{p_k\}_{k=1}^\infty$ a sequence of projections 
with $p_{k}\in \al_{n_k},\; k\in{\mathbb N}$. Suppose that
\begin{align}\label{cons}
\lim_{k\to\infty}\lV
(1-p_k)f_1(H_{1,n_k})\cdots f_m(H_{m,n_k})
\rV=\lim_{k\to\infty}\lV
f_1(H_{1,n_k})\cdots f_m(H_{m,n_k})(1-p_k)
\rV=0,
\end{align}
for all continuous functions $f_1,\cdots,f_m$
over $\mathbb R$ with ${\rm supp}f_i\subset U_i$,
$i=1,\cdots,m$.
Then we have
\begin{align*}
\liminf_{k}\frac 1{2n_k+1} \log {\rm Tr}_{n_k} p_k\ge
\sup\left\{
\mu(x)\; :\; x\in \prod_{i=1}^m U_i
\right\}.
\end{align*}
\end{lem}
{\it Proof} 
First we show $\lim_{k\to\infty}\omega(p_k)=1$, for any ergodic state
$\omega$ over $\al$ with $\omega(A_i)\in U_i,\; i=1,\cdots , m$.
Let $f_1,\cdots,f_m$ be continuous functions 
over $\mathbb R$ 
with
$0\le f_i\le 1$, $f_i(\omega(A_i))=1$, 
and ${\rm supp}f_i\subset U_i$.
From von Neumann's ergodic Theorem
\cite{simon}, by the ergodicity of
$\omega$, we have
\begin{align*}
\lim_{n\to\infty} \omega\lmk
f_1(H_{1,n})\cdots f_m(H_{m,n})
\rmk=\prod_{i=1}^mf_i(\omega(A_i))=1.
\end{align*}
From this and the assumption (\ref{cons}),
 we have
\begin{align*}
1=
\liminf_{k\to\infty}\lv
\omega\lmk
f_1(H_{1,n_k})\cdots f_m(H_{m,n_k})
\rmk\rv\le
\liminf_{k\to\infty}\lv
\omega(p_kf_1(H_{1,n_k})\cdots f_m(H_{m,n_k})p_k)\rv
\le\liminf_{k\to\infty}\omega(p_k).
\end{align*}
We thus obtain $\lim_{k\to\infty}\omega(p_k)=1$.
\\
Let $\omega$ be an ergodic state with $\omega(A_i)\in U_i,\; i=1,\cdots , m$.
The above assertion means
for any $0<\varepsilon<1$,
$\omega(p_k)\ge 1-\varepsilon$ for $k$ large enough.
We thus have $\beta_{\varepsilon,n_k}(\omega)\le \log {\rm Tr}_{n_k}p_k$
eventually.
Applying Theorem \ref{Bj},
we have
\begin{align}\label{lbe}
s(\omega)=\liminf_{n\to\infty}\frac{1}{2n+1}
\beta_{\varepsilon,n}(\omega)\le
\liminf_k\frac{1}{2n_k+1}\log {\rm Tr}_{n_k} p_k.
\end{align}
Now we claim that this inequality can be extended to
general translation invariant states. 
To do so, we use a standard technique
in statistical mechanics.(See \cite{simon},\cite{BR96}):
Let $\omega$ be a translation invariant state with
$\omega(A_i)\in U_i,\; i=1,\cdots , m$.
Define a translation invariant state $\bar\omega_L$, by
\begin{align*}
\bar\omega_L:=\frac{1}{2L+1}\sum_{j\in[-L,L]}
\tilde\omega_L\circ\gamma_j,\quad
\tilde\omega_L:=\bigotimes\omega\vert_{\al_L},
\end{align*}
for each $L\in{\mathbb N}$.
It is well known that $\bar\omega_L$ is ergodic and
$\bar\omega_L\to \omega$ in weak$*$-topology
of $S_{\gamma}(\al)$.
In particular, $\bar\omega_L(A_i)\in U_i,\;i=1,\cdots,m$ eventually,
as $L\to\infty$.
Furthermore, the mean entropy $s(\bar\omega_L)$ of
$\bar\omega_L$ is equal to $\frac{1}{2L+1}S_{[-L,L]}(\omega)$,
where $S_{[-L,L]}(\omega)$ is the von Neumann entropy
of $\omega\vert_{\al_L}$.\\
Applying (\ref{lbe}) for the ergodic state
$\bar\omega_L$, we get
\begin{align*}
\frac{1}{2L+1}S_{[-L,L]}(\omega)=s(\bar\omega_L)
\le
\liminf_k\frac{1}{2n_k+1}\log {\rm Tr}_{n_k} p_k.
\end{align*}
Taking $L\to\infty$ limit, we obtain
\begin{align*}
s(\omega)\le\liminf_k\frac{1}{2n_k+1}\log {\rm Tr}_{n_k} p_k.
\end{align*}
This implies the result. $\square$
%By
%\begin{align*}
%\sup\{s(\omega)\; :\; \omega\in S_\gamma(\al),\;
%\omega(A_i)\in U_i,\; i=1,\cdots, m\}\\=
%\sup\{\sup\{
%s(\omega)\; :\; \omega\in S_\gamma(\al),\;
%\omega(A_i)=y_i
%\}\; : \; y=(y_1,\cdots,y_m)\in\prod_{i=1}^m U_i,\; i=1,\cdots, m\}\\
%=
%\sup\{\mu(y)\; : \;  y\in \prod_{i=1}^m U_i\},\end{align*}
%we obtain the result.$\square$\\
\section{An ideal tower in $\prod_{k}M_{n_k}({\mathbb C})/
\oplus M_{n_k}({\mathbb C})$}\label{tower}
%\section{An ideal $I_s$ of 
%$\prod_{k}M_{n_k}({\mathbb C})/
%\oplus_k M_{n_k}({\mathbb C})$
%}\label{is}
Let $\{n_k\}_k$ be a subsequence of $\mathbb N$.
We fix this sequence in the rest of this section.
Define a $C^*$-algebra $B$ by
\begin{align*}
B:=\prod_{k}M_{n_k}({\mathbb C})
=
\left\{
\lmk x_k\rmk\; :\;
\sup_{k}
\lV x_k\rV<\infty,\;
x_k\in M_{n_k},\; k\in{\mathbb N}
\right\},
\end{align*}
and its closed ideal $D$ by
\begin{align*}
D:=\oplus_k M_{n_k}({\mathbb C})=\left\{
\lmk x_k\rmk\in B\; :\;
\lim_{k\to\infty}
\lV x_k\rV=0
\right\}.
\end{align*}
For each $i=1,\cdots, m$,
$(H_{i,n_k})$ is a self-adjoint element in $B$.
We denote the quotient map from $B$ to $A:=B/D$ by $\pi$.
The $C^*$-algebra $A$ has real rank zero.
It is well known that for any projection $p$ in $A$, 
there exists a projection
 $(p_k)$ in $B$
such that $\pi((p_k))=p$.
Similarly, for any partial isometry $v$ in $A$,
there exists a partial isometry $(v_k)$
in $B$ such that $\pi((v_k))=v$. (See Theorem 1.3 \cite{lin2}.)\\
Take $R>\max_{i=1,\cdots,m}\lV A_i\rV$ and define a compact subset
$X$ of $\br^m$ by $X:=\prod_{i=1}^m[-R,R]$.
Note that $X_s=\{x\in X\;;\; \mu(x)\ge s\}$.
As the sequences $H_{i,n},\;i=1,\cdots,m$ mutually commute 
asymptotically, i.e.,
\begin{align*}
\lim_{n\to\infty}
\lV
\left[H_{i,n},H_{j,n}\right]\rV=0,\quad i,j=1,\cdots, m,
\end{align*}
we have 
\begin{align*}
\left[\pi((H_{i,n_k})),\pi((H_{j,n_k}))\right]=0,\quad
i,j=1,\cdots,m.
\end{align*}
Therefore, we can define
a $*$-homomorphism $\varphi:C(X)\to A$ by
\begin{align*}
\varphi(f):=f\lmk
\pi((H_{1,n_k})),\cdots,\pi((H_{m,n_k}))
\rmk,\quad f\in C(X).
\end{align*}
For each $i=1,\cdots,m$, define 
 $h_i\in C(X)$ by
\[
h_i(x)=x_i,\quad x=(x_1,\cdots,x_m).
\]
Clearly, $\varphi(h_i)=\pi((H_{i,n_k}))$.\\\quad\\
We define a closed subset $S$ of $X$ as follows :
$x\in S$ iff 
for any neighborhood $U$ of $x$, there exists
$f\in C(X)$ with ${\rm supp }f\subset U$ such that $\varphi(f)\neq 0$.
\begin{lem}
\[
{\rm dom}\mu =S.
\]
\end{lem}
{\it Proof}
To prove $S\subset {\rm dom}\mu$,
let $x\in S$.
%Then for any $\varepsilon>0$, there exists 
%$f\in C(X)$ with
% $\varphi(f)\neq0$ and ${\rm supp}f\subset 
%B_{\varepsilon(x)}$.\\
Fix $\varepsilon>0$. 
We first prove that
there exists a translation invariant state 
$\bar\omega_\varepsilon$
such that 
$(\bar\omega_\varepsilon(A_1),\cdots, 
\bar\omega_\varepsilon(A_m))\in B_{4\varepsilon\sqrt m}(x)$.
By Lemma \ref{epl}, there exists a projection
$p$ in $A$ such that 
\[
\hat\varphi(1_{\overline{B_{\varepsilon(x)}}})
\le p\le \hat\varphi(1_{B_{2\varepsilon(x)}}),
\]
where $\hat\varphi:C(X)^{**}\to{\al}^{**}$ is the extension of $\varphi$.
Let $(p_k)$ be a projection in $B$
such that $\pi((p_k))=p$.\\
From the definition of $p$,
we have $\lV \varphi(h_i)p-x_ip\rV<3\varepsilon$.
This implies
\begin{align}\label{acc}
\lV \lmk H_{i,n_k}-x_{i}\rmk p_{k}\rV<3
\varepsilon,
\end{align}
for $k$ large enough.
By the assumption $x\in S$, there exists
$f\in C(X)$ with
 $\varphi(f)\neq0$ and ${\rm supp}f\subset 
B_{\varepsilon}(x)$. We may assume $0\le f\le 1$.
As $0\le f\le1_{\overline{B_{\varepsilon(x)}}} $,
we have 
$0\le\varphi(f)\le \hat\varphi(1_{\overline{B_{\varepsilon(x)}}})
\le p$.
From this, there exists a positive element $(a_{k})$ in $B$
such that $0\le a_{k}\le p_{k} \;\forall k$,
and $\pi((a_{k}))=\varphi(f)$.
As $\pi((a_{k}))=\varphi(f)\neq 0$, there exists
a subsequence $(a_{k_M})_M$ of $(a_k)$ such that 
$a_{k_M}\neq 0$, for all $M$.
Let $\omega_M$ be a state over 
$\al_{n_{k_M}}$ with a density matrix
$
\frac{a_{k_M} }{{\rm Tr}_{n_{k_M}}a_{k_M} }
$,
and define a translation invariant state 
$\bar \omega_M$ by
\[
\bar \omega_M:=\frac{1}{2n_{k_M}+1}
\sum_{j\in [-n_{k_M},n_{k_M}]}\tilde\omega_M\circ\gamma_j,\quad
\tilde\omega_M:=\bigotimes
\omega_M.
\]
Then we have 
\[
\bar \omega_M(A_i)
=\omega_M(H_{i,n_{k_M}}),
\] for all $i=1,\cdots,m$.
On the other hand, from (\ref{acc}),
we get 
\begin{align*}
\lv
\omega_M(H_{i,n_{k_M}})-x_{i}
\rv
=\lv\frac{{\rm Tr}_{n_{k_M}}\lmk a_{k_M} p_{k_M}
\lmk H_{i,n_{k_M}}-x_{i}\rmk\rmk}{{\rm Tr}_{n_{k_M}}a_{k_M} }
\rv\le
\lV
p_{k_M}
\lmk H_{i,n_{k_M}}-x_{i}\rmk
\rV<4\varepsilon,
\end{align*}for $M$ large enough.
Hence, 
we obtain
\begin{align*}
\lv
\bar \omega_M(A_i)-x_i\rv<4\varepsilon,
\end{align*}for $M$ large enough.
We define $\bar\omega_{\varepsilon}:=\bar\omega_M$,
for such large $M$, and the claim is established.\\
Next, we consider the net of translation
invariant states $\{\bar\omega_{\varepsilon}\}_{\varepsilon>0}$ 
taken as above.
As the space of translation invariant states
$S_\gamma(\al)$ is wk$*$-compact,
the net $\{\bar\omega_{\varepsilon}\}_{\varepsilon}$ 
has a convergent subnet 
$\{\bar\omega_{\varepsilon'}\}_{\varepsilon'},\;
\bar\omega_{\varepsilon'}\to \omega\in S_\gamma(\al)$,
in wk$*$-topology.
As 
$(\bar\omega_\varepsilon(A_1),\cdots, 
\bar\omega_\varepsilon(A_m))\in B_{\delta}(x)$
eventually for any $\delta>0$, we have
\[
\omega(A_i)=x_i,\; i=1,\cdots,m.
\] 
Hence we obtain a translation invariant state $\omega$
with $\omega(A_i)=x_i,\; i=1,\cdots,m$.
This implies $x\in{\rm dom }\mu$.
We thus obtain $S\subset {\rm dom }\mu$.\\
Next we prove ${\rm dom}\mu\subset S$.
Let $x\in S^c$. By the definition of $S$,
there exists an open neighborhood $U$ of $x$ such that
$\varphi(f)=0$ for all $f\in C(X)$ with ${\rm supp}f\subset U$.
We claim $\mu(x)=-\infty$. Assume $\mu(x)\in{\br}$.
Then, there exists $\omega\in S_{\gamma}(\al)$
such that $\omega(A_i)=x_i,\; i=1,\cdots,m$.
Define $\rho:=\bigotimes\omega\vert_{\al_{0}}$.
This state is an ergodic state with $\rho(A_i)=\omega(A_i)
=x_i,\; i=1,\cdots,m$.
Choose $\varepsilon>0$ and $f_i\in C([-R,R]),\; i=1,\cdots,m$,
so that
$\overline{B_{\sqrt m \varepsilon}(x)}\subset U$,
and $f_i(x_i)=1,\; f_i\vert_{B_{\varepsilon}(x_i)^c}=0,\; 0\le f_i\le 1$.
Then $f\in C(X)$ given by
\[
f(y):=f_1(y_1)\cdots f_m(y_m),\quad
y=(y_1,\cdots,y_m)\in \prod[-R,R]=X
\]
has its support in $\overline{B_{\sqrt m \varepsilon}(x)}$,
hence in $U$.
Therefore, by the assumption
we have 
\[0=
\varphi(f)=f_1(\pi((H_{1,n_k})))\cdots f_m(\pi((H_{m,n_k})))
=\pi\lmk
\lmk f_1(H_{1,n_k})\cdots f_m(H_{m,n_k})
\rmk\rmk.
\]
This means
\[
\lim_{k\to\infty}f_1(H_{1,n_k})\cdots f_m(H_{m,n_k})=0.
\]
On the other hand, by the ergodicity of 
$\rho$, we get
\[
0=\lim_k \rho\lmk f_1(H_{1,n_k})\cdots f_m(H_{m,n_k})\rmk
=f_1(\rho(A_1))\cdots f_m(\rho(A_m))=
f_1(x_1)\cdots f_m(x_m)=1,
\]
which is a contradiction. Hence we obtain $\mu(x)=-\infty$.
$\square$\\\\
%From this Lemma, we have
%\[
%X_s:=\{x\in \br^m\;:\; \mu(x)\ge s\}=\{x\in X\;:\; \mu(x)\ge s\}
%\]
%for all $s\in\br^m$.\\
Given $C^*$-algebras $A_1$, $A_2$ and
a $*$-homomorphism $\rho : A_1\to A_2$,
we extend $\rho$ naturally to a $*$-homomorphism
from $M_N(A_1)$ to $M_N(A_2)$ for each $N\in{\mathbb N}$,
and denote it by the same symbol $\rho$.
\begin{prop}\label{lem19}
For each $s\in{\br^m}$, define a set of projections 
${\cal S}_s$ in $A$ by
\[
{\cal S}_s:=\left\{
e\in {\rm Proj} A\; :\;
\exists (e_k)\in {\rm Proj}(B)\; s.t.\; e=\pi((e_k)),
\limsup\frac{1}{2n_k+1}\log {\rm Tr}_{n_k} e_k<s
\right\}
\]
and let $I_s$ be the closed ideal
of $A$ generated by ${\cal S}_s$.
Then the following statements hold:
\begin{description}
\item[(i)]
For any $p\in {\rm Proj}(A)$,
$p$ is in the ideal $I_s$ iff 
\[
\limsup_k\frac{1}{2n_k+1}\log {\rm Tr}_{n_k} p_k<s,
\]for all $(p_k)\in {\rm Proj} B$
such that $p=\pi((p_k))$.
\item[(ii)]
For all $N\in{\mathbb N}$ and $p\in {\rm Proj}(M_N(I_s))$,
\[
\limsup_k\frac{1}{2n_k+1}\log {\rm Tr}_{n_k} p_k<s,
\]for all $(p_k)\in {\rm Proj} M_N(B)$
such that $p=\pi(p_k)$.
\item[(iii)]
Let $p$ be a projection in $A$, and suppose that
there exists $(p_k)\in{\rm Proj}(B)$ such that
$p=\pi((p_k))$ and
\begin{align*}
\liminf_k\frac{1}{2n_k+1}\log {\rm Tr}_{n_k} p_k\ge s.
\end{align*} 
Then for any $N\in{\mathbb N}$ and 
$q\in {\rm Proj}(M_N(I_s))$,
we have\[
q\lesssim p.
\]
\item[(iv)]
For any $x\in X_s$ and its open neighborhood $U$ in $X$,
there exists a continuous function $g\in C(X)$
with $g\vert_{U^c}=0,\; 0\le g\le 1$, satisfying the following
property:
for any $p\in {\rm Proj}(A)$ such that $0\le\varphi(g)\le p$ and
$(p_k)\in {\rm Proj}(B)$ with $p=\pi((p_k))$,
we have
\[
\liminf \frac{1}{2n_k+1}\log {\rm Tr}_{n_k} p_k\ge s.
\]
\item[(v)]
If $g\in C(X)$ satisfies $g\vert_{X_s}=0$, then
$\varphi(g)\in I_s$.
\item[(vi)]
If $s<0$, then $I_s=\{0\}$.
\end{description}
\end{prop}
{\it Proof}
(i)The "If" part is trivial.
To prove the "only if" part, 
we first note, as shown in \cite{lin2}, that $e\in{\cal S}_s$
iff $\limsup\frac{1}{2n_k+1}\log {\rm Tr}_{n_k} e_k<s$
for {\it any} $(e_k)\in {\rm Proj}(B)$ such that $e=\pi((e_k))$.
This follows from the fact that
for any $(e_k), (f_k)\in{\rm Proj}B$ with $\pi((e_k))=\pi((f_k))$,
$e_k\sim f_k$ holds eventually as $k\to\infty$.
From Proposition 1.13 of \cite{lin2},
for any $p\in {\rm Proj}(I_s)$,
there exist finite number of projections $e_1,\cdots,e_l$
in ${\cal S}_s$ such that $p\lesssim e_1\oplus\cdots\oplus e_l$
in $M_l(I_s)$. Let $v\in M_{1,l}(I_s)$ be a partial isometry 
such that 
\[
vv^*=p,\quad v^*v\le e_1\oplus\cdots\oplus e_l.
\]
For this $v$, there exist $(p_k)\in {\rm Proj}(B)$,
$(q_k)\in {\rm Proj}(M_l(B))$, and $(v_k)\in M_{1,l}(B)$
such that $p_k=v_kv_k^*, v_k^*v_k\le q_k$,
$\pi(p_k)=p,\; \pi(q_k)= e_1\oplus\cdots\oplus e_l$,
and $\pi(v_k)=v$.
This can be proven by the same argument as Theorem 1.3 of
\cite{lin2}.
As $e_i\in {\cal S}_s,\; i=1,\cdots, l$,
there exists $(e_{k}^i)\in {\rm Proj}(B)$ such that
$\pi((e_k^i))=e_i$ with $\limsup_k\frac{1}{2n_k+1}\log {\rm Tr}_{n_k} e_k^i< s$.
For these $(e_k^i)$, we have
$\pi(e^1_k\oplus\cdots\oplus e^l_k)
=e_1\oplus\cdots\oplus e_l=\pi(q_k)$.
Therefore, for $k$ large enough, we have
\[
{\rm Tr}_{n_k} p_k={\rm Tr}_{n_k} v_k^*v_k\le {\rm Tr}_{n_k} q_k=\sum_{i=1}^l{\rm Tr}_{n_k}(e_k^i).
\]
From this, we obtain
\begin{align*}
\limsup_k\frac{1}{2n_k+1}\log {\rm Tr}_{n_k} p_k
\le
\limsup_k\frac{1}{2n_k+1}\log \sum_{i=1}^l{\rm Tr}_{n_k} e_k^i
=\max_{1\le i\le l}\limsup_k\frac{1}{2n_k+1}\log {\rm Tr}_{n_k} e_k^i
<s.
\end{align*}
This means $p\in{\cal S}_s$.
From the assertion at the beginning of the proof, we obtain
the claim.
\\
(ii)Using the fact that $M_N(I_s)$ is the closed ideal of $M_N(A)$ generated
by ${\cal S}_s^N:=
\{(e_1\oplus\cdots\oplus e_N)\; :\; e_i\in{\cal S}_s \}$,
 proof of (ii) is the same as that of (i).\\
(iii)Let $p$ and $q$ be as in (iii). From (ii), for any projection 
$(q_k)$ in $M_N(B)$ with $q=\pi((q_k))$, we have
\[
s_0:=\limsup_k\frac{1}{2n_k+1}\log {\rm Tr}_{n_k}q_k<s.
\]
With this and the assumption on $p$, we get
\begin{align*}
\liminf_k\frac{1}{2n_k+1}\log \frac{{\rm Tr}_{n_k} p_k}{{\rm Tr}_{n_k} q_k}
\ge \liminf_k\frac{1}{2n_k+1}\log {{\rm Tr}_{n_k} p_k}-
\limsup_k\frac{1}{2n_k+1}\log{{\rm Tr}_{n_k} q_k}
\ge s-s_0>0.
\end{align*}
This means ${\rm Tr}_{n_k}p_k>{\rm Tr}_{n_k}q_k$, hence
$q_k\lesssim p_k$ for $k$ large enough.
Therefore, we have $q\lesssim p$.\\
(iv) 
For $x\in{X_s}\subset{\rm dom}\mu$ and its open neighborhood $U$ in $X$,
there exist open subsets $U_i,\; i=1,\cdots,m$ of $[-R,R]$
with $x\in\overline{\prod_{i=1}^m U_i}\subset U$,
and $g\in C(X)$ satisfying $0\le g\le 1$,
$g\vert_{\overline{\prod_{i=1}^m U_i}}=1$,
$g\vert_{U^c}=0$.
We prove that this $g$ enjoys the required property.
%Let $p$ be a projection
%such that $0\le \varphi(p)\le p$.
For any $f_i\in C(\br)$ with 
$supp f_i\subset U_i$, we have 
$(1-g)\prod_{i=1}^mf_i\circ h_i=0$
and thus 
$(1-\varphi(g))\prod_{i=1}^m
f_i(\pi((H_{i,n_k})))=(1-\varphi(g))\prod_{i=1}^m
\varphi(f_i\circ h_i)=0.$
\\
Now, for any $p\in {\rm Proj}(A)$ with $0\le 
\varphi(g)\le p$,
we get 
\[
(1-p)\prod_{i=1}^m
f_i(\pi((H_{i,n_k})))=0,
\]
from 
\begin{align*}
0\le 
\lmk
\pi\lmk 
f_1(H_{1,n_k})\cdots f_m(H_{m,n_k})\rmk\rmk^*
(1-p)
\lmk
\pi\lmk 
f_1(H_{1,n_k})\cdots f_m(H_{m,n_k})\rmk\rmk\\
\le \lmk
\pi\lmk 
f_1(H_{1,n_k})\cdots f_m(H_{m,n_k})\rmk\rmk^*
(1-\varphi(g))
\lmk
\pi\lmk 
f_1(H_{1,n_k})\cdots f_m(H_{m,n_k})\rmk\rmk=0.
\end{align*}
For any projection $(p_k)$ in 
$B$ with $\pi((p_k))=p$,
%this implies
%\[
%\pi((1-p_k)
%f_1(H_{1,n_k})\cdots f_m(H_{m,n_k}))=0.
%\]
%%and $\pi(
%%f_1(H_{1,n_k})\cdots f_m(H_{m,n_k})(1-p_k))=0$, 
this means
\[
\lim_k\lV(1-p_k)
f_1(H_{1,n_k})\cdots f_m(H_{m,n_k})\rV=0.\]
Similarly, we have \[\lim_k\lV f_1(H_{1,n_k})\cdots 
f_m(H_{m,n_k})(1-p_k)\rV=0.
\]
Applying Lemma \ref{lem10}, we obtain
\begin{align*}
\liminf\frac{1}{2n_k+1}\log {\rm Tr}_{n_k} p_k
\ge\sup\{\mu(y)\; :\; y\in \prod_{i=1}^m U_i\}
\ge \mu(x)\ge s.
\end{align*}
(v)
Let $x$ be an element in $X_s^c\cap X$,
and $\varepsilon>0$ a positive number
such that $\overline{{B_{\sqrt m\varepsilon}(x)}}\subset X_s^c$.
Let $g\in C(X)$ be a function $0\le g\le 1$ 
with ${\rm supp}g\subset B_{\varepsilon}(x)$.
We prove $\varphi(g)\in I_s$.
It suffices to consider the case $\varphi(g)\neq 0$.
\\
From Lemma \ref{epl}, there exists a projection $r$ in $A$
such that $\varphi(g)\le\hat\varphi(1_{{\rm supp}g})\le
 r\le \hat{\varphi}(1_{B_{\varepsilon}(x)})$.
For this $r$, we have 
\begin{align}\label{rv}
\lV\pi(((H_{i,n_k}-x^i)r_k))\rV=\lV\varphi(h_i)r-x^ir\rV\le \varepsilon,
\end{align}
where $(r_k)$ is a projection in $B$ such that $r=\pi((r_k))$.
As $\varphi(g)\neq 0$, we have $r\neq 0$.
Therefore, there exists a subsequence $\{r_{k_M}\}$
of $\{r_k\}$ consisting of all the nonzero projections in $\{r_k\}$.
For this subsequence, and
for any $\delta>0$,
\[
\lmk
\frac{{\rm Tr}_{n_{k_M}} r_{k_M} H_{1,n_{k_M}}}{{\rm Tr}_{n_{k_M}} r_{k_M}},\cdots,
\frac{{\rm Tr}_{n_{k_M}} r_{k_M} H_{m,n_{k_M}}}{{\rm Tr}_{n_{k_M}} r_{k_M}}
\rmk\in \overline{ \lmk B_{\sqrt m\varepsilon}(x)\rmk}_\delta,
\]eventually from (\ref{rv}).
Therefore from Lemma \ref{lem9}, we obtain
\begin{align*}
\limsup_k\frac 1{2n_k+1} \log {\rm Tr}_{n_k} r_{n_k}\le
\sup\{\mu(x):x\in \overline{ B_{\sqrt m \varepsilon}(x)}\}.
\end{align*}
By the upper semi-continuity of $\mu$, we have
\[
\sup\{
\mu(x)\; : \; x\in \overline{ B_{\sqrt m \varepsilon}(x)}
\}<s.
\]
This means $r\in I_s$. As $\varphi(g)\le r$,
we have $\varphi(g)\in I_s$.\\
General cases follow from this, using partition of unity
and approximation of $g$ with continuous functions with
supports in $X_s^c$.\\
(vi)
Assume $s<0$.
If $I_s\neq \{0\}$, then there exists a nonzero projection $e\in {\cal S}_s$.
Let $(e_k)\in {\rm Proj}(B)$
such that $e=\pi(e_k)$ and
\[
\limsup_k\frac{1}{2n_k+1}\log {\rm Tr}_{n_k} e_k<s
.\]
As $e\neq 0$,  
there exists a subsequence $(e_k')$ of $(e_k)$
such that ${\rm Tr}_{n_{k'}}e_k'\ge 1$.
Therefore, we have
\[
s<0\le \limsup_k\frac{1}{2n_k+1}\log {\rm Tr}_{n_k} e_k<s,
\]
which is a contradiction. Therefore, $I_s=\{0\}$
$\square$\\\quad\\
Now we construct an ideal tower.
\begin{prop}\label{prop19}
Let $\eta>0$ be a positive number. Then
\begin{description}
\item[(i)]There exists a finite sequence of
real numbers $s_0>s_1>\cdots>s_n$,
such that
\begin{align*}
s_0=\sup\left\{
\mu(x)\;:\; x\in{\br}^m
\right\}<\infty,\quad s_n<0,\\
X_{s_k}\subset B_\eta(X_{s_{k-1}}),\quad k=1,\cdots,n,\\
{\rm dom}\mu= X_0= X_{s_n}.
\end{align*}
Furthermore, $X_{s_k},\; k=1,\cdots, n$ are compact and convex,
and $X_{s_0}$ consists of one point
\[
x_0:=(\frac{{\rm Tr}_{0} A_1}{d},\cdots \frac{{\rm Tr}_{0} A_m}{d})\in\br^m.
\]
\item[(ii)]For $s_0,\cdots,s_n$ in (i),
there exist points $\lambda_{ij}\in X,\; i=0,\cdots,n-1,\;j=1.\cdots,l_i$
with $\lambda_{ij}\neq\lambda_{i'j'}$ for $(i,j)\neq(i',j')$,
such that
\[
\lambda_{ij}\in X_{s_i}\backslash X_{s_{i-1}},\;
i=0,1,\cdots, n-1,\; j=1,\cdots,l_i,
\]
where we set $X_{s_{-1}}:=\phi$. For
each $k=0,\cdots,n-1$, the set $\{\lambda_{ij}\;:\;
i=0,\cdots ,k,\;j=1,\cdots,l_i\}$
is $2\eta$-dense in $X_{s_{k+1}}$.
Furthermore, $l_0=1$ and $X_{s_0}=\{\lambda_{01}\}$.
\item[(iii)]
For $\{\lambda_{ij}\}$ in (ii) and any $\beta>0$,
there exist mutually orthogonal projections $\{r_{ij}\; : i=0,
\cdots,n-1,j=1,\cdots,l_i\}$ in $A$ with
$r_{ij}\in I_{s_{i-1}},\quad i=1,\cdots, n-1$,
$1-r_{01}\in I_{s_0}$, satisfying the following conditions: for any $g\in C(X)$,
\begin{align}\label{gest}
&r_{ij}\varphi(g)r_{i'j'}=0,\quad (i,j)\neq(i',j'),\nonumber\\
&\lV
\varphi(g)r_{ij}-g(\lambda_{ij})r_{ij}
\rV\le
\sup\left\{
\lv g(\zeta)-g(\lambda_{ij})\rv\; :\;
\lv \zeta-\lambda_{ij}\rv<\beta\;, \zeta\in X
\right\},
\end{align}
and for a projection $r:=\sum_{ij}r_{ij} $,
we have
\begin{align}\label{gest2}
\lV
\varphi(g)r-\sum_{ij}g(\lambda_{ij})r_{ij}
\rV\le
\sup\left\{
\lv g(\zeta)-g(\lambda)\rv\; :\;
\lv \zeta-\lambda\rv<\beta,\; \zeta,\lambda\in X
\right\}.
\end{align}
Furthermore, for each $i=0,\cdots, n-1, j=1,\cdots,l_i$,
there exists a projection $(r_{ij}^k)$ in $B$
with $\pi((r_{ij}^k))=r_{ij}$, such that
\begin{align}\label{rijk}
\liminf\frac{1}{2n_k+1}\log {\rm Tr}_{n_k} r_{ij}^k\ge s_i.
\end{align}
\end{description}
\end{prop}
{\it Proof} (i) is proven in Lemma \ref{lem21}.
To prove (ii), choose for each $k=0,\cdots,n-1$,
 a finite set of elements 
$E_k:=\{\zeta_j^{(k)}\}_{j=1,\cdots,m_k}$ in $X_{s_k}$
which is $\eta$-dense in the compact set $X_{s_{k}}$.
As $X_{s_{k+1}}\subset B_\eta(X_{s_{k}})$,
$E_k$ is $2\eta$-dense in $X_{s_{k+1}}$, $k=0,\cdots,n-1$.
Define $\Lambda_i:=\lmk
\cup_{k=0}^{n-1}E_k\rmk \cap(X_{s_i}\backslash X_{s_{i-1}})$,
$i=0,\cdots,n-1$.
Then we have $\Lambda_i\subset X_{s_i}\backslash X_{s_{i-1}}$
and $\cup_{i=0}^k\Lambda_i$
is $2\eta$-dense in $X_{s_{k+1}}$.
Labeling elements in $\Lambda_i$ as $\Lambda_i=\{\lambda_{ij}\}_{j=1,\cdots,l_i}$,
for each $i=0,\cdots, n-1$, we obtain
$\{\lambda_{ij}\}$ which satisfy the conditions in (ii).\\
Now for an arbitrary $\beta>0$, we construct projections
$\{r_{ij}\}$ in (iii).
Fix $\delta>0$ so that
\[
\delta<\frac 14\min\{\lv
\lambda_{ij}-\lambda_{i'j'}\rv\; :\;
(ij)\neq (i'j')\}\wedge
\frac 14\min\{\inf_{\lambda\in X_{s_{i-1}}}
\lv \lambda_{ij}-\lambda \rv)\; :\;
i=1,\cdots, n-1\}
\wedge\frac 14\eta
\wedge\frac 14\beta.
\]
For each $i=0,\cdots,n-1,\;j=1,\cdots,l_i$,
by Lemma \ref{epl}, there exists a projection $r_{ij}$
in $A$ such that 
\begin{align}\label{rij}
\hat\varphi(1_{\overline{B_\delta(\lambda_{ij})}})
\le r_{ij}\le \hat\varphi(1_{B_{2\delta}(\lambda_{ij})}).
\end{align}
As $B_{2\delta}(\lambda_{ij})\cap
B_{2\delta}(\lambda_{i'j'})=\phi$ for $(ij)\neq(i'j')$,
 these inequalities imply
that $\{r_{ij}\}$ are mutually orthogonal, and
 (\ref{gest}), (\ref{gest2}) hold.\\
To see $r_{ij}\in I_{s_{i-1}},\; i=1,\cdots,n-1$,
let $g\in C(X)$ be a function such that $0\le
g\le 1,\; g\vert_{\overline{B_{2\delta}(\lambda_{ij})}}=1,\;
g\vert_{B_{3\delta}(\lambda_{ij})^c}=0$.
As $\delta$ is taken small enough so that $X_{s_{i-1}}\subset
B_{3\delta}(\lambda_{ij})^c$, we have $g\vert_{X_{s_{i-1}}}=0$.
From Proposition \ref{lem19} (v), this implies 
$\varphi(g)\in I_{s_{i-1}}$.
By (\ref{rij}), we have
\[
r_{ij}\le \hat\varphi(1_{B_{2\delta}(\lambda_{ij})})\le
\varphi(g).
\] Hence we obtain $r_{ij}\in I_{s_{i-1}}$.
To see $1-r_{01}\in I_{s_0}$, define
$g\in C(X)$ to be a function such that $0\le
g\le 1,\; g\vert_{\overline{B_{\frac{\delta}{2}}
(\lambda_{01})}}=1,\;
g\vert_{B_{\delta}(\lambda_{01})^c}=0$.
Then we have $(1-g)\vert_{X_{s_0}}=0$.
Therefore, by (v) of Proposition {\ref{lem19}}, we obtain $\varphi(1-g)\in I_{s_0}$.
This and the inequality 
$\varphi(g)\le \hat\varphi
(1_{\overline{B_{\delta}(\lambda_{01})}})\le r_{01}$
implies $1-r_{01}\in I_{s_0}$.\\
To show (\ref{rijk}), we apply (iv) of Proposition \ref{lem19}.
to $\lambda_{ij}\in X_{s_i}$, $B_{\delta}(\lambda_{ij})$, and
 obtain 
$g\in C(X)$ such that $g\vert_{B_{\delta}(\lambda_{ij})^c}=0,\; 0\le g\le 1$.
As $0\le \varphi(g)\le \hat\varphi
(1_{\overline{B_{\delta}(\lambda_{ij}})})\le r_{ij}$,
there exists 
a projection $(r_{ij}^k)$ in $B$
such that $r_{ij}=\pi((r_{ij}^k))$ and
$\liminf\frac{1}{2n_k+1}\log Tr_{n_k} r_{ij}^k\ge s_i$.
$\square$
\section{Proof of Theorem \ref{main}}\label{proof}
In this section, we prove Theorem \ref{main}.
\begin{defn}\label{DF}
Let $X$ be a compact metric space.
For a finite subset $\cal F$ of $C(X)$, we say that
$X$ satisfies the condition $D_{\cal F}$ if for any $\varepsilon>0$,
there exist a positive number $\delta:=\delta_D(\varepsilon,{\cal F},X)>0$
and a positive integer $N:=N_{\cal D}(\varepsilon,{\cal F},X)$
satisfying the following:
For any unital $C^*$-algebra $\cal B$, unital $*$-homomorphism
$\varphi :C(X)\to {\cal B}$ and a projection $p\in {\cal B}$
satisfying
\[
\lV p\varphi(f)-\varphi(f)p\rV<\delta,\quad \forall f\in
{\cal F},
\]
there exist $m,r\in{\mathbb N}$,
$\xi_j\in X,\; j=1,\cdots, m$,
$\lambda_l\in X,\; l=1,\cdots, r$,
two sets of mutually orthogonal projections 
$p_j,\in {\rm Proj}(M_N(p{\cal B}p))\; j=1,\cdots, m$
and $q_l,\in {\rm Proj}(M_{N+1}(p{\cal B}p))\;l=1,\cdots, r$
with
\[
\sum_{j=1}^m p_j=1_{M_N(p{\cal B}p)},\quad
\sum_{l=1}^r q_l=1_{M_{N+1}(p{\cal B}p)}
\]
such that
\[
\lV
p\varphi(f)p\oplus\sum_{j=1}^m f(\xi_j)p_j
-\sum_{l=1}^r f(\lambda_l)q_l
\rV<\varepsilon,
\]for all $f\in{\cal F}$.
\end{defn}
%In particular, this holds when $X$ consists of one point.
\begin{thm}[\cite{eglp}]
Consider the compact metric space
\[
I^n=[-1,1]\times [-1,1]\times\cdots\times
[-1,1],
\] and 
let ${\cal F}:=\{g_1,\cdots,g_m\}$
be a set of generators of $C(I^n)$.
Then $I^n$ satisfies condition ${\cal D}_{\cal F}$.
\end{thm}
%\begin{rem}
%The proof in \cite{eglp} applies for the case that $X$ consists of
%one point.
%\end{rem}
Now, a nonempty compact convex subset in ${\br}^m$ is
homeomorphic to $I^l$ for some $0\le l\le m$.
Therefore, each $X_{s_{k}}\;1\le k\le n$ 
in Proposition \ref{prop19}
satisfies the condition 
${\cal D}_{{\cal F}\vert_{X_{s_k}}}$
for ${\cal F}=\{1,h_1,\cdots,h_m\}$.\\
The following Lemma can be proven following the idea 
of \cite{gl}. We give a sketch of its proof in
Appendix.
\begin{lem}\label{lem25}
Let $X$ be a compact metric space and $X_0,\cdots,X_n$
a finite sequence of its closed subsets such that
\[
\{x_0\}=X_0\subset X_1\subset\cdots \subset X_n\subset X,
\]
where $x_0$ is an element in $X$.
Let $\cal F$ be a finite subset of $C(X)$, and assume that
each $X_k,\; k=1,\cdots,n$ satisfies the condition
$D_{{\cal F}_k}$ for ${\cal F}_k:=
\{f\vert_{X_k}\; f\in {\cal F}\}\subset\subset
C(X_k)$.
Furthermore, let $A$ be a unital $C^*$-algebra with real rank zero
and
$I_0,\cdots, I_n$ a finite sequence of its closed ideals with
\[
\{0\}=I_n\subset\cdots \subset I_1\subset I_0\subset A,
\]
where $I_{k+1}$ is an ideal of $I_k$ for $k=0,\cdots,n-1$.
Let $\pi_k:A\to A/I_{k}$
and $\pi_{k,k+1}: A/{I_{k+1}}\to A/I_k$ be the quotient maps.\\
Suppose that there exists a unital $*$-homomorphism 
$\varphi:C(X)\to A$ satisfying
\begin{align}\label{vsp}
\pi_k\circ\varphi(g)=0,\quad {\rm if}\quad g\vert_{X_k}=0,\;
g\in C(X),
\end{align}
for each $k=1,\cdots, n$.
Then for any $\varepsilon>0$, there exists a positive number
$\delta=\delta_{T}(\varepsilon,{\cal F},\{X_k\}_{k=0}^n)>0$
satisfying the following: if $p$ is a projection in $A$
with
\begin{align}\label{sct}
\pi_0(p)=1,\quad \lV[\varphi(f),p]\rV<\delta,\quad\forall f\in{\cal F},
\end{align} then
there exist a sequence of positive integers
$N_0,\cdots,N_{n-1}\in{\mathbb N}$, 
$*$-homomorphisms 
$h_k\; :\; C(X)\to M_{N_k}(I_k)\; k=0,\cdots, n-1$,
and $H\; : \; C(X)\to M_{N_0+\cdots+N_{n-1}+1}(A)$ 
with finite dimensional range,
such that
\begin{align}\label{do}
\lV\bar p\varphi(f)\bar p\oplus h_0(f)\oplus
\cdots \oplus h_{n-1}(f)-H(f)\rV<\varepsilon,\quad
\forall f\in{\cal F}.
\end{align}
Furthermore, $h_k$ and $H$ are of the form
\begin{align*}
&h_k(f)=\sum_{j=1}^{L_k}f(\xi_{kj})p_{kj},\\
&p_{kj}\in {\rm Proj}(M_{N_k}(I_k))\quad
j=1\cdots L_k,
\quad{\rm mutually\;orthogonal},\\
&\xi_{kj}\in X_{k+1},\quad
L_k\in{\mathbb N},\; k=0,\cdots, n-1,
\end{align*}
and
\begin{align*}
&H(f)=\sum_{i=1}^Lf(\zeta_i)q_i,\\
&q_i\in {\rm Proj}(M_{N_0+\cdots+N_{n-1}+1}(A)),
\quad i=1,\cdots
L,\quad {\rm mutually\;orthogonal},\\
&\zeta_i\in X,\quad L\in{\mathbb N},
\end{align*}
with
\begin{align}\label{quni}
\bar p\oplus\sum_{j=1}^{L_0} p_{0j}\oplus
\cdots \oplus \sum_{j=1}^{L_{n-1}} 
p_{{n-1},j}=\sum_{i=1}^L q_i.
\end{align}
Here, we used the notation $\bar p:=1-p$.
\end{lem}
Combining all the results so far obtained,
we can show that $\varphi$ can be approximated
by a $*$-homomorphism with finite dimensional range:
%in the following sense:
\begin{thm}\label{finiteapp}
For any $\varepsilon>0$, there exists
a unital $*$-homomorphism $G:C(X)\to A$ with finite
dimensional range 
such that
\[
\lV \varphi(f)-G(f)\rV
<\varepsilon,\quad {\rm for\; all}\; f\in {\cal F}=\{1,h_1,\cdots,h_m\}.
\]
\end{thm}
{\it Proof} 
We follow the argument in \cite{gl}.
Fix $\varepsilon>0$. As $X$
is compact, there exists $\eta>0$
such that
$\lv f(\zeta)-f(\lambda)\rv<\frac{\varepsilon}{8}$
for all $f\in{\cal F}$ and $\zeta,\lambda\in X$
with $\lv \zeta-\lambda\rv<2\eta$.\\
For this $\eta>0$, we can find a
finite sequence $s_0>s_1>\cdots>s_n$
of real numbers and
$\lambda_{ij}\in X$ satisfying the conditions 
in (i) and (ii) of Proposition \ref{prop19}.
Put $X_k:=X_{s_k},\; k=0,\cdots,n$.
\\
Each $X_k=X_{s_{k}},\;1\le k\le n$ 
in Proposition \ref{prop19}
satisfies the condition 
${\cal D}_{{\cal F}\vert_{X_{s_k}}}={\cal D}_{{\cal F}\vert_{X_k}}$
for ${\cal F}$.
For the unital $C^*$-algebra $A$, we obtain
an ideal tower
$\{0\}=I_{n}\subset\cdots\subset I_{1}\subset I_{0}\subset A$
where $I_k:=I_{s_k},\; k=0,\cdots, n$. Note that
$I_{{k+1}}$ is an ideal of $I_{k}$.
By $s_n<0$, we have $I_n=I_{s_n}=\{0\}$ from Proposition \ref{lem19} (vi).
Let $\pi_k:A\to A/I_{k},\; \pi_{k,k+1}:A/I_{{k+1}}\to A/I_{k},\quad
k=0,\cdots,n-1$ 
be the quotient maps.
By Proposition \ref{lem19} (v), we have
$\pi_k\circ \varphi(g)=0$, for all $g\in C(X)$
with $g\vert_{X_{k}}=0$, for each $k\ge 1$.
Applying Lemma \ref{lem25}, we obtain a positive number
$\delta_{T}(\frac {\varepsilon}{3},{\cal F},\{X_k\}_{k=0}^n)$.\\
As $X$ is compact, there exists $\beta>0$
such that $\lv f(\zeta)-f(\lambda)\rv<\frac 13
\delta_{T}(\frac {\varepsilon}{3},{\cal F}, \{X_k\}_{k=0}^n)\wedge
\frac\varepsilon{50}$, for all $f\in{\cal F}$ and
$\zeta,\lambda\in X$ with $\lv \zeta-\lambda\rv<\beta$.\\
For this $\beta$, applying Proposition \ref{prop19},
we can find mutually orthogonal
projections $r_{ij}\in {\rm Proj}(A)$ satisfying conditions in 
(iii) of Proppsition \ref{prop19}.
By the choice of $\beta$, we have
%\begin{align*}
%\lV
%\varphi(f)r_{ij}-f(\lambda_{ij})r_{ij}
%\rV\le
%\frac 13
%\delta_{T}(\frac {\varepsilon}{3},{\cal F}, \{X_k\}_{k=0}^n)\wedge
%\frac\varepsilon{50}\le
%\frac\varepsilon{50},\quad\forall f\in{\cal F},
%\end{align*}
%and
\begin{align}\label{rloc}
\lV
\varphi(f)r-\sum_{ij}f(\lambda_{ij})r_{ij}
\rV\le
\frac 13
\delta_{T}(\frac {\varepsilon}{3},{\cal F}, \{X_k\}_{k=0}^n)\wedge
\frac\varepsilon{50}\le
\frac 13
\delta_{T}(\frac {\varepsilon}{3},{\cal F}, \{X_k\}_{k=0}^n),\quad\forall f\in{\cal F},
\end{align}
where we put $r:=\sum_{ij}r_{ij}$.
From this inequality, we get
\begin{align}\label{rvc}
\lV
\lcm
\varphi(f),r
\rcm
\rV<
\frac 23 \delta_{T}(\frac {\varepsilon}{3},{\cal F}, \{X_k\}_{k=0}^n)\wedge
\frac{2\varepsilon}{50}
<
\delta_{T}(\frac {\varepsilon}{3},{\cal F}, \{X_k\}_{k=0}^n),\quad\forall f\in{\cal F}.
\end{align}
Furthermore, we have $\pi_0(r)=1$.
Applying Lemma \ref{lem25},
we obtain a sequence of positive integers
$N_0,\cdots,N_{n-1}\in{\mathbb N}$, 
$*$-homomorphisms 
$h_k\; :\; C(X)\to M_{N_k}(I_k),\; k=0,\cdots, n-1$,
and $H\; : \; C(X)\to M_{N_0+\cdots+N_{n-1}+1}(A)$
such that
\begin{align}\label{do1}
\lV\bar r\varphi(f)\bar r\oplus h_0(f)\oplus
\cdots \oplus h_{n-1}(f)-H(f)\rV<\frac{\varepsilon}{3},\quad
\forall f\in{\cal F},
\end{align}
with $\bar r=1-r$.
Furthermore, $h_k$ and $H$ are of the form
\begin{align*}
&h_k(f)=\sum_{l=1}^{L_k}f(\xi_{kl})p_{kl},\\
&p_{kl}\in {\rm Proj}(M_{N_k}(I_k))\quad
l=1\cdots L_k,
\quad{\rm mutually\;orthogonal}\\
&\xi_{kl}\in X_{k+1},\quad
L_k\in{\mathbb N},\; k=0,\cdots, n-1
\end{align*}
and
\begin{align*}
&H(f)=\sum_{i=1}^Lf(\zeta_i)q_i\\
&q_i\in {\rm Proj}(M_{N_0+\cdots+N_{n-1}+1}(A)),\quad i=1,\cdots L,\quad {\rm mutually\;orthogonal}\\
&\zeta_i\in X,\quad L\in{\mathbb N},
\end{align*}
where the projections satisfy
\begin{align}\label{quni1}
\bar r\oplus\sum_{j=1}^{L_0} p_{0j}\oplus
\cdots \oplus \sum_{j=1}^{L_{n-1}} p_{n-1,j}
=\sum_{i=1}^L q_i.
\end{align}
Now recall that for each $k=0,\cdots,n-1$,
the set
$\{\lambda_{ij}\}_{i=0,\cdots,k,j=1,\cdots,l_i}$ is $2\eta$-dense in
$X_{{k+1}}$.
Therefore, for each $\xi_{kl}\in X_{{k+1}}$,
we can find $\lambda_{i'j'}$ with $0\le  i'\le k$ such 
that \[\lv
\lambda_{i'j'}-\xi_{kl}\rv<2\eta.
\]
By the choice of $\eta$, this means
\[
\max\{\lv f(\lambda_{i'j'})-f(\xi_{kl})\rv\;:\; f\in{\cal F}\}
<\frac\varepsilon 8.
\]
Choose such $\lambda_{i'j'}$ for each $\xi_{kl}$
and denote it by $\hat\lambda(\xi_{kl})$.
Let $\hat q_{ij}^k\in {\rm Proj}(M_{N_k}(A))$, 
$k=0,\cdots,n-1,i=0,\cdots,n-1,j=1,\cdots,l_i$ be projections given by
\begin{align*}
\hat q_{ij}^k:=\ltm
\begin{gathered}
\sum_{l\;:\; \hat\lambda(\xi_{kl})=\lambda_{ij}}p_{kl},\quad
i\le k\\
0,\quad i\ge k+1
\end{gathered}.\right.
\end{align*}
As each $p_{kl}$ is in $M_{N_k}(I_{k})$,
$\hat q_{ij}^k$ is in $M_{N_k}(I_k)$ as well. 
Note that 
\[
\sum_{i=0,\cdots,n-1,j=1,\cdots,l_i}
\hat q_{ij}^k
=\sum_{i=0}^k \sum_{j=1}^{l_i}\hat q_{ij}^k
=\sum_{i=0}^k\sum_{j=1}^{l_i}
\sum_{l\;:\; \hat\lambda(\xi_{kl})=\lambda_{ij}}p_{kl}
=\sum_{l=1}^{L_k}p_{kl}.
\]
For each $k=0,\cdots ,n-1$, we define a
$*$-homomorphism $h'_k\;:\; C(X)\to M_{N_k}(I_k)$ by
\begin{align*}
h'_k(g):=\sum_{i=0}^k\sum_{j=1}^{l_i}g(\lambda_{ij})
\lmk
\sum_{l\;:\; \hat\lambda(\xi_{kl})=\lambda_{ij}}p_{kl}\rmk
=\sum_{i=0}^{n-1}\sum_{j=1}^{l_i}g(\lambda_{ij})\hat q_{ij}^k,
\quad g\in C(X).
\end{align*}
From the choice of $\hat\lambda(\xi_{kl})$,
we have
\begin{align}\label{hhd}
\lV
h_k(f)-h_k'(f)
\rV
%=\lV
%\sum_{i=0}^k
%\sum_{j=1}^{l_i}\sum_{l\;:\; \hat\lambda(\xi_{kl})=\lambda_{ij}}
%\lmk
%f(\xi_{kl})-f(\hat\lambda(\xi_{kl}))
%\rmk p_{kl}
%\rV\nonumber \\
=\max\ltm
\lv
f(\xi_{kl})-f(\hat\lambda(\xi_{kl}))
\rv\; :\; l=1,\cdots,L_k
\rtm<\frac\varepsilon 8,\quad \forall f\in{\cal F}.
\end{align}
Define mutually orthogonal projections
$Q_{ij},\; i=0,\cdots ,n-1,\; j=1,\cdots l_i$
in $M_{N_0+\cdots+N_{n-1}}(A)$ by
\begin{align}
Q_{ij}:=\hat q_{ij}^0\oplus\cdots\oplus
\hat q^i_{ij}\oplus \hat q^{i+1}_{ij}\oplus
\cdots \hat q^{n-1}_{ij}\nonumber\\
=0\oplus\cdots\oplus 0\oplus
\hat q^i_{ij}\oplus \hat q^{i+1}_{ij}\oplus
\cdots \hat q^{n-1}_{ij}.
\end{align}
As $\hat q^k_{ij}$ is in $M_{N_k}(I_k)$, 
$Q_{ij}$ is in the ideal $M_{N_0+\cdots+N_{n-1}}(I_i)$.
For this $Q_{ij}$ have
\begin{align}\label{hdd}
\sum_{i=0}^{n-1}\sum_{j=1}^{l_i}g(\lambda_{ij})Q_{ij}
%=\sum_{i=0}^{n-1}\sum_{j=1}^{l_i}g(\lambda_{ij})
%\lmk
%\hat q_{ij}^1\oplus\cdots\oplus
%\hat q^i_{ij}\oplus \hat q^{i+1}_{ij}\oplus
%\cdots \hat q^{n-1}_{ij}\rmk\nonumber\\
=h_0'(g)\oplus\cdots\oplus h_{n-1}'(g),
\end{align}
for all $g\in C(X)$.
Furthermore, we have
\begin{align*}
\bar r\oplus\sum_{ij}Q_{ij}
=\bar r\oplus \sum_l p_{0l}\oplus\cdots\oplus
\cdots \sum_l p_{n-1,l}
=\sum_i q_i=H(1).
\end{align*}
From ({\ref {do1}}), (\ref{hhd}), and (\ref{hdd}),
we obtain
\begin{align}\label{rhh}
\lV
\bar r \varphi(f)\bar r \oplus 
\sum_{i=0}^{n-1}\sum_{j=1}^{l_i}f(\lambda_{ij})Q_{ij}-H(f)
\rV
%=\lV
%\bar r \varphi(f)\bar r \oplus h_0'(f)\oplus\cdots
%\oplus h_{n-1}'(f)
%-H(f)
%\rV\nonumber\\
%\le
%\lV
%0\oplus \lmk h_0(f)-h_0'(f)\rmk
%\oplus\cdots\oplus \lmk h_{n-1}(f)-h_{n-1}'(f)\rmk
%\rV\nonumber\\
%+\lV
%\bar r \varphi(f)\bar r \oplus h_0(f)\oplus\cdots\oplus h_{n-1}(f)
%-H(f)
%\rV\nonumber\\
%<\frac{\varepsilon}{8}+\frac{\varepsilon}{3}
<\frac{2\varepsilon}{3},\quad\forall f\in {\cal F}.
\end{align}
Now, recall that for each $r_{ij}$,
there exists $(r_{ij}^k)\in {\rm Proj}B$ with 
$r_{ij}=\pi((r_{ij}^k))$, which satisfies
\[
\liminf_k\frac{1}{2n_k+1}\log {\rm Tr}_{n_k}r_{ij}^k
\ge s_i.\]
Then by (iii) of Proposition \ref{lem19}, 
for any $N\in{\mathbb N}$ and 
any $q\in {\rm Proj}(M_N(I_{i}))$,
we have\[
q\lesssim r_{ij}.
\]
In particular, we have
\[
Q_{ij}\lesssim r_{ij}.
\]
This means there exists a partial
isometry $v_{ij}\in M_{1,N_0+\cdots+N_{n-1}}(A)$
such that 
\begin{align}\label{vr}
r_{ij}':=v_{ij}v_{ij}^*\le r_{ij},\quad
Q_{ij}=v_{ij}^*v_{ij}.
\end{align}
Let $v$ be a partial isometry given by
\[
v:=\lmk \bar r,\sum_{ij}v_{ij} \rmk
\in M_{1,N_0+\cdots+N_{n-1}+1}(A),
\]
%By the relation (\ref{vr}),
%we get
%\begin{align}\label{vv}
%v^*v=\bar r \oplus \sum_{ij}Q_{ij}=H(1),\quad
%vv^*=\bar r+\sum_{ij}r_{ij}'.
%\end{align}
%Furthermore,
%we have
%\begin{align}\label{mog}
%v\lmk
%\bar r \varphi(f)\bar r \oplus 
%\sum_{i=0}^{n-1}\sum_{j=1}^{l_i}f(\lambda_{ij})Q_{ij}\rmk v^*
%=\bar r \varphi(f)\bar r + 
%\sum_{i=0}^{n-1}\sum_{j=1}^{l_i}f(\lambda_{ij})r_{ij}'.
%\end{align}
and define $G:C(X)\to A$ by
\[
G(f):=vH(f)v^*+\sum_{ij}f(\lambda_{ij})(r_{ij}-r_{ij}').
\]
%By (\ref{vv}), $vH(f)v^*$
%and $\sum_{ij}g(\lambda_{ij})(r_{ij}-r_{ij}'),\quad f,g\in C(X)$ are orthogonal
%and
%\[
%vH(f)v^*vH(g)v^*=vH(f)H(1)H(g)v^*=vH(fg)v^*.
%\]
It is easy to check that 
$G$ is a unital $*$-homomorphism with finite dimensional
range.\\
%Furthermore, by
%\[
%G(1)=vv^*+\sum_{ij}(r_{ij}-r_{ij}')
%=\bar r+\sum_{ij}r_{ij}'+\sum_{ij}(r_{ij}-r_{ij}')=1,
%\]
%$G$ is unital.\\
%Now 
%we estimate $\lV \varphi(f)-G(f)\rV$, for $f\in{\cal F}$.
By (\ref{rloc}), (\ref{rvc}) and (\ref{rhh}),
%\begin{align*}
%\lV \varphi(f)-G(f)\rV\\
%\le \lV \varphi(f)-\lmk
%r\varphi(f)r+\bar r\varphi(f)\bar r
%\rmk\rV+
%\lV\lmk
%r\varphi(f)r+\bar r\varphi(f)\bar r
%\rmk
%-G(f)\rV\\
%\le
%2\lV\lcm \varphi(f),r\rcm\rV
%+\lV
%r\lmk
%\varphi(f)r-\sum_{ij}f(\lambda_{ij})r_{ij}
%\rmk
%\rV\\
%+\lV
%\bar r\varphi(f)\bar r+\sum_{ij}f(\lambda_{ij})r_{ij}-G(f)
%\rV
%<\frac 2{25}\varepsilon+\frac{\varepsilon}{50}
%+\lV
%\bar r\varphi(f)\bar r+\sum_{ij}f(\lambda_{ij})r_{ij}
%-G(f)\rV,\quad\forall f\in{\cal F}.
%\end{align*}
%Now by (\ref{mog}),
%\begin{align*}
%\lV
%\bar r\varphi(f)\bar r+\sum_{ij}f(\lambda_{ij})r_{ij}
%-G(f)\rV
%=\lV
%\bar r\varphi(f)\bar r
%-vH(f)v^*
%+\sum_{ij}f(\lambda_{ij})r_{ij}'\rV\\
%=\lV v\lmk
%\bar r\varphi(f)\bar r
%\oplus \sum_{ij}f(\lambda_{ij})Q_{ij}-H(f)\rmk v^*\rV
%<\frac{2\varepsilon}{3},\quad\forall f\in{\cal F}.
%\end{align*}
we obtain
\begin{align*}
\lV \varphi(f)-G(f)\rV<\varepsilon.\quad
\forall f\in{\cal F}.\quad\square
\end{align*}
{\it Proof of Theorem \ref{main}}\\
The proof is by contradiction. Assume the assertion were false.
Then there exists $\varepsilon>0$ and a subsequence $\{n_k\}_k$
of $\mathbb N$ such that 
\begin{align}\label{cont}
\inf\ltm
\max_{1\le i\le m}\lV
H_{i,{n_k}}-X_{i,{n_k}}
\rV\;:\; \lcm X_{i,{n_k}}, X_{j,{n_k}}
\rcm=0,\;X_{i,{n_k}}\in (\al_{n_k})_{sa},\;
i,j=1,\cdots,m
\rtm\ge \varepsilon.
\end{align}
Applying Theorem \ref{finiteapp} to this subsequence
$\{n_k\}_k$ and $A=\prod_{k}M_{n_k}({\mathbb C})/
\oplus M_{n_k}({\mathbb C})$,
we obtain
a unital $*$-homomorphism $G:C(X)\to A,$ with
 finite dimensional range
such that
\[
\lV \varphi(f)-G(f)\rV<\frac{\varepsilon}{2},\quad
\forall f\in{\cal F}=\{1,h_1,\cdots,h_m\}.
\]
We can represent $G$ as
\[
G(g)
=\sum_{j=1}^Mg(\zeta_j)Q_j,\quad g\in C(X),
\]where $\zeta_j\in X$, and mutually orthogonal projections
$Q_j$ in $A$.
From Lemma \ref{lem23}, there exist mutually orthogonal projections
$(Q^k_j)$ in $B$
such that $Q_j=\pi((Q_j^k))$,
and $\sum_jQ_j^k=1_{\al_{n_k}}$.
By $\varphi(h_i)=\pi((H_{i,n_k}))$,
we get
\begin{align*}
\lV
\pi\lmk
(H_{i,n_k}-\sum_{j=1}^Mh_i(\zeta_j)Q_j^k)
\rmk
\rV
=\lV
\varphi(h_i)-G(h_i)
\rV<\frac{\varepsilon}{2},\quad i=1,\cdots,m.
\end{align*}
This means, for $k$ large enough,
we have
\begin{align*}
\max_{1\le i\le m}\lV
H_{i,n_k}-\sum_{j=1}^Mh_i(\zeta_j)Q_j^k
\rV
<\frac{\varepsilon}{2},\quad i=1,\cdots,m.
\end{align*}
Define $X_{i,k}$ by
\[
X_{i,k}:=\sum_{j=1}^Mh_i(\zeta_j)Q_j^k\in(\al_{n_k})_{sa}.
\]
Then we have sequences $(X_{1,k})_k,\cdots,(X_{m,k})_k$ such that
\[
\lcm
X_{i,k},X_{j,k}
\rcm=0,\quad X_{i,k}\in \lmk
\al_{n_k}\rmk_{sa},
\]
and
\[
\max_{1\le i\le m}\lV
H_{i,{n_k}}-X_{i,k}
\rV<\frac{\varepsilon}{2},
\]
eventually.
This contradicts (\ref{cont}). $\square$\\\\
\noindent
{\bf Acknowledgement.}\\
{The author is grateful for Professor Hiroki Matui for helpful 
explanation on their
papers \cite{matui} and \cite{ms}  
and the result of G.
Gong and H. Lin \cite{gl}.
The present research is 
supported by JSPS 
Grant-in-Aid for Young Scientists (B), 
Sumitomo Foundation, and Inoue Science Research Award.}
\appendix
\section{$C^*$-algebra of real rank zero}
In this section, we list the results on $C^*$-algebra 
of real rank zero
that we use.
\begin{lem}\label{epl}
Let $X$ be a compact metric space,
$A$ a unital $C^*$-algebra with real rank zero,
and $\varphi:C(X)\to A$ a unital $*$-homomorphism.
Then for any closed subset $V$ of $X$ 
and open subset $U$ of $X$
with $V\subset U$, there exists a projection $r$ in $A$
such that
\[
\hat\varphi(1_{V})\le r\le \hat\varphi(1_U).
\]
Here, $\hat \varphi$ is the homomorphism from $C(X)^{**}$
to $A^{**}$ given as the unique extension of
$\varphi$.
\end{lem}
{\it Proof} See \cite{bowen}.
%Let $p:=\hat \varphi (1_V), q:=1-\hat\varphi(1_U)\in A^{**}$.
%Let $B$ be a hereditary $C^*$-subalgebra of $A$ 
%generated by $\varphi(C_0(U))$, i.e., 
%$B=\overline{\varphi(C_0(U))A\varphi(C_0(U))}$.
%There exists a function $f\in C(X)$ such that $f\vert_V=1,f\vert_{U^c}=0,\; 0\le f\le 1$.
%The element $a:=\varphi(f)$ is in $B$.
%We have $p\le a:=\varphi(f)\le 1-q$.
%As $A$ has real rank zero, there exists $s\in{\rm Proj}(B)$
%such that $\lV a(1-s)\rV<1$.\\
%Let $x:=(1-a)^{\frac 12}(1-s)$. 
%Then we have $x^*x\in (1-s)A(1-s)$, and
%\[
%1-s-x^*x=1-s-(1-s)(1-a)(1-s)=(1-s)a(1-s)\le
%\lV
%(1-s)a(1-s)
%\rV(1-s), 
%\]
%with $\lV(1-s)a(1-s)\rV<1$.
%Hence $\lv x\rv=(x^*x)^{\frac 12}$ is an invertible element in $(1-s)A(1-s)$,
%i.e., there exists $y\in (1-s)A(1-s)$ such that
%$\lv x\rv y=y\lv x\rv=1-s$.
%Then for $u:=xy$, we have
%\[
%x=x(1-s)=x y\lv x\rv=u\lv x\rv,
%\]
%and
%\[
%u^*u=y^*x^*xy=1-s,\quad
%uu^*uu^*=uu^*.
%\]
%Hence we have a polar decomposition
%$x=u\lv x\rv$, where $u$ is a partial isometry of $A$
%and $u^*u=1-s$.\\
%For $r:=1-uu^*$, we have
%\[
%1-r=uu^*=xyy^*x^*\le
%\lV y\rV^2xx^*\le
%\lV y\rV^2(1-a)\le \lV y\rV^2(1-p).
%\] 
%From this, we get
%\[
%1-r\le 1-p.
%\]
%As $q\le 1-a$, we have 
%\[
%q\le q(1-a)q\le q,
%\]
%which implies
%\[
%qaq=0.
%\]
%This implies 
%\[
%a^{\frac 12}q=0.
%\]
%From this, we get
%\[
%1-r= uu^*\ge u\lv x\rv^2u^*=xx^*=(1-a)^{\frac 12}(1-s)(1-a)^{\frac 12}
%\ge (1-a)^{\frac 12}q(1-a)^{\frac 12}=q.
%\]
%Hence we obtain a projection $r\in{\rm Proj}(A)$
%such that $p\le r\le 1-q$.
$\square$
\begin{lem}\label{lem20}
Let $A$ be a unital $C^*$-algebra with real rank zero,
and $I$ a closed ideal of $A$ with quotient map $\pi : A\to A/I$.
Let $h$ be a positive element in $A$ with $\pi(h)^2=\pi(h)$,
and $B$ a hereditary $C^*$-subalgebra of $A$ generated by $h$.
Then there exists $p\in{\rm Proj}(B)$ such that
$\pi(h)=\pi(p)$.
\end{lem}
{\it Proof} See \cite{zhang}.
$\square$
\begin{lem}\label{lem23}
Let $A$ be a unital $C^*$-algebra with real rank zero.
Let $I$ be a closed ideal of $A$ and
$\pi\;: \; A\to A/I$ the quotient map.
Then for any mutually orthogonal projections
 $\{p_l\}_{l=1}^N$ in $A/I$, and a projection $p$
in $A$ with $p_l\le \pi(p),\; l=1,\cdots, N$,
there exist mutually orthogonal projections
$\tilde p_l\in A ,\; l=1,\cdots, N$ 
such that $\pi(\tilde p_l)=p_l$
and $\tilde p_l\le p$.
Furthermore, if $\pi(p)=\sum_lp_l$,
$\tilde p_l$s can be taken to satisfy $\sum_l \tilde p_l=p$.\end{lem}
{\it Proof} See \cite{lin1}.
$\square$
\begin{lem}\label{lem24}
Let $A$ be a unital $C^*$-algebra with real rank zero.
Let $I$ be a closed ideal of $A$ with quotient map
$\pi:A\to A/I$.
Let $\{p_l\}_{l=1}^N$ be mutually orthogonal projections
in $A$, and put $p:=\sum_{l=1}^Np_l$.
Then for any $\delta>0$ and $x_1,\cdots,x_m\in pAp$ 
with
\begin{align}\label{delta}
\lV \pi(x_i)\rV<\delta,\quad i=1,\cdots, m,
\end{align}
there exist $e\in {\rm Proj}(I)$ and $e_l\in {\rm Proj}(p_lIp_l)$
such that 
\[
e=\sum_{l=1}^N e_l,\quad\lV
x_je-ex_j
\rV<4\delta,
\quad \lV
(1-e)x_j(1-e)
\rV<2\delta,\quad j=1,\cdots,m.
\]
\end{lem}
{\it Proof} See \cite{gl}. $\square$
%From (\ref{delta}), there exist $a_j\in pIp,\; j=1,\cdots, m$
%such that $\lV x_j-a_j\rV<\delta$.
%For each $1\le l\le N$, $p_l Ip_l$
%is a hereditary $C^*$-subalgebra of $A$.
%As $A$ has real rank zero, for $p_la_ja_j^*p_l,p_la_j^*a_jp_l\in p_lIp_l$,
%$j=1,\cdots,m$, 
%there exists $e_l\in{\rm Proj}(p_lIp_l)$
%such that
%\begin{align*}
%\lV
%e_lp_l a_ja_j^*p_l-p_l a_ja_j^*p_l\rV<
%\lmk
%\frac{\delta}{2N^2}
%\rmk^2,\\
%\lV
%p_l a_j^*a_jp_l e_l-p_l a_j^*a_jp_l\rV<
%\lmk
%\frac{\delta}{2N^2}
%\rmk^2.
%\end{align*}
%Then for $e:=\sum_{l=1}^N e_l\in{\rm Proj}(I)$, we have
%\begin{align*}
%\lV
%a_j-ea_je
%\rV=\lV
%pa_jp-ea_je
%\rV
%\le
%\sum_{l,l'=1}^N
%\lV
%(1-e_l)p_la_jp_{l'}
%\rV
%+\lV
%e_la_jp_{l'}(1-e_{l'})
%\rV\\
%\le
%\sum_{l,l'=1}^N\lV
%(1-e_l)p_la_ja_j^*p_l(1-e_l)
%\rV^{\frac 12}
%+\lV
%(1-e_l')p_{l'}a_j^*a_jp_{l'}(1-e_{l'})
%\rV^{\frac 12}
%<\delta,
%\end{align*}
%for all $j=1,\cdots,m$.
%From this we obtain
%\begin{align*}
%\lV x_je-ex_j\rV
%\le
%\lV(x_j-a_j)e-e(x_j-a_j)\rV
%+\lV a_je-ea_j\rV\\
%\le
%2\lV x_j-a_j\rV
%+\lV
%(1-e)a_je
%\rV
%+\lV
%ea_j(1-e)
%\rV<2\cdot\delta+2\delta=4\delta,
%\end{align*}
%and
%\begin{align*}
%\lV
%(1-e)x_j(1-e)
%\rV\le
%\lV
%(1-e)(x_j-a_j)(1-e)
%\rV+\lV
%(1-e)a_j(1-e)
%\rV
%<\delta+\delta=2\delta.\quad \square
%\end{align*}
%
%
%
%
\section{Proof of Lemma \ref{pmu}}
The proof follows the standard arguments in statistical
mechanics, relating the mean entropy and the free energy.
(See\cite{BR96}).\\
First we prove 
$p(\alpha)\ge
\mu(x)+(\alpha,x)$
for all $x\in\br^m$ and $\alpha\in\br^m$.
This is trivial if $\mu(x)=-\infty$.
If  $\mu(x)>-\infty$, then
for any $\varepsilon>0$, there exists
a state $\omega\in S_\gamma(\al)$ satisfying
\begin{align}\label{smue}
\omega(A_i)=x_i,\;1\le i\le m, \quad 
s(\omega)\ge \mu(x)-\varepsilon.
\end{align}
Using the positivity of the relative entropy 
$0\le S(\omega\vert_{\al_n},
\frac{e^{(2n+1)\sum_{i=1}^m\alpha_i H_{i,n}}}
{{\rm Tr}_{n}e^{(2n+1)\sum_{i=1}^m\alpha_i H_{i,n}}})
$, we get
\[
\frac{1}{2n+1}S_{[-n,n]}(\omega)
+\sum_{i=1}^m\alpha_i \omega(A_i)
\le
p(\alpha),\quad\forall n\in{\mathbb N}.
\]
Taking $n\to\infty$ limit, we obtain
\[
s(\omega)+\sum_{i=1}^m\alpha_i \omega(A_i)
\le p(\alpha).
\]
By (\ref{smue}),  
we get
$p(\alpha)\ge s(\omega)+\sum_{i=1}^m\alpha_ix_i\ge
\mu(x)+(\alpha,x)-\varepsilon\;,x\in{\br}^m,\;\alpha\in{\br}^m$.
We thus obtain
\[
p(\alpha)\ge 
\sup\{\mu(x)+(\alpha,x)\;:\;x\in{\br}^m\}.
\]
\\
Next we prove $p(\alpha)\le
\sup\left\{
\mu(x)+(\alpha,x)\;:x\in{\br}^m
\right\}$ for all $\alpha\in{\br}^m$.
We define a state $\rho$ on $\md$ by
\begin{align*}
\rho:=\frac{e^{\sum_{i=1}^m\alpha_i A_i}}
{{\rm Tr}e^{\sum_{i=1}^m\alpha_i A_i}}.
\end{align*}
From this, we can define a translation invariant
state $\tilde\rho:=\bigotimes_{\mathbb Z}\rho$.
We can easily see that for 
$y:=(\tilde\rho( A_1),\cdots\tilde\rho( A_m))$,
\begin{align*}
p(\alpha)=s(\tilde\rho)
+\sum_{i=1}^m\alpha_i\tilde\rho( A_i)
\le
\mu(y)+(\alpha,y)
\le
\sup\left\{
\mu(x)+(\alpha,x)\;:x\in{\br}^m
\right\}.
\end{align*}
Hence we obtain the first equality.\\
To prove the second assertion, we recall the following fact:
For a function $G : \br^m\mapsto [-\infty,\infty]$, we define
its Legendre transform $G^*:\br^m\to[-\infty,\infty]$ by
\[
G^{*}(u):=\sup\{
(\alpha,u)-G(\alpha)\; :\; \alpha\in{\br}^m
\},\quad u\in\br^m.
\]
\begin{thm}\cite{convex}
Let $F$ be a convex and lower semi-continuous function of $\br^m$
into $(-\infty,\infty]$.
Then we have
\[
F=F^{**}.
\]
\end{thm}
Applying this theorem to $F:=-\mu$,
we obtain the claim.
$\square$
\section{Proof of Lemma \ref{lem25}}
{\it Proof}
It suffices to show the claim for the case
that $\lV f\rV\le 1$, for all $f\in{\cal F}$.
For a fixed $\varepsilon>0$,
we define a finite sequence of positive numbers
$\delta_0,\cdots,\delta_n$ inductively by
$\delta_n:=\varepsilon$ and $\delta_{k}:=\min\{
\frac{1}{100}\delta_D(\frac{\delta_{k+1}}{3},{\cal F}_{k+1},X_{k+1}),
\frac 1{100}\delta_{k+1},\frac 1{10}\},k=0,\cdots, n-1$.
We take $\delta<\frac 12 \delta_0$.
Let $p$ be a projection satisfying (\ref{sct}) for this $\delta$.
We consider the following proposition 
$(A_k),\; k=0,\cdots, n-1$:\\
\begin{center}
$(A_k)$ : There exist positive integers 
$N_{j}^k\in{\mathbb N},\; j=0,\cdots, k$ and
$*$-homomorphisms with finite dimensional range
$h_j^k:C(X)\to M_{N_j^k}(I_j/I_{k+1})$ for
$j=0,\cdots, k$, and
$H^k\; : \; C(X)\to M_{N_0^k+\cdots+N_k^{k}+1}
(A/I_{k+1})$ satisfying
\begin{align}\label{pif}
\lV\pi_{k+1}(\bar p\varphi(f)\bar p)
\oplus h_0^k(f)\oplus
\cdots \oplus h_k^k(f)-H^k(f)\rV<\delta_{k+1},\quad
\forall f\in{\cal F}.
\end{align}
Furthermore, $h^k_j$ and $H^k$ are of the form
\begin{align*}
&h_j^k(f)=\sum_{i=1}^{L^k_j}f(\xi_{ji}^k)p_{ji}^k,\\
&p_{ji}^k\in {\rm Proj}(M_{N^k_j}(I_j/I_{k+1}))\quad
i=1,\cdots, L^k_j,\; 
\quad{\rm mutually\;orthogonal}\\
&\xi_{ji}^k\in X_{j+1},\quad j=0,\cdots,k,
\end{align*}
and
\begin{align*}
&H^k(f)=\sum_{i=1}^{L^k}f(\zeta_i^k)q_i^k\\
&q_i^k\in {\rm Proj}(M_{N_0^k+\cdots+N_k^k+1}(A/I_{k+1})),\quad i=1,\cdots
L^k,\quad {\rm mutually\;orthogonal}\\
&\zeta_i^k\in X,\quad L^k\in{\mathbb N},
\end{align*}
with
\begin{align*}
\pi_{k+1}(\bar p)\oplus\sum_{i=1}^{L_0^k}p_{0i}^k
\oplus\cdots
\oplus 
\sum_{i=1}^{L^k_k}p_{ki}^k
=\sum_{i=1}^{L^k} q_i^k.
\end{align*}
\end{center}
$(A_{n-1})$ corresponds to the claim of the Lemma.
\\
We assume that $(A_k)$ holds and prove 
that $(A_{k+1})$ holds, for $k=0,\cdots,n-2$.
Applying Lemma \ref{lem23},
to $\{p_{ji}^k\}_{i=1,\cdots, L^k_j}\subset
M_{N_j^k}(I_j/I_{k+1})$ in $(A_k)$, we can find 
mutually orthogonal projections
$\{\tilde p_{ji}^k\}_{i=1,\cdots, L^k_j}$
in $M_{N_j^k}(I_j/I_{k+2})$ such that
$\pi_{k+1,k+2}(\tilde p_{ji}^k)=p^k_{ji}$.
%As we have $\pi_{j,k+2}(\tilde p_{ji}^k)=\pi_{j,k+1}\circ
%\pi_{k+1,k+2}(\tilde p_{ji}^k)=0$,
%$\tilde p_{ji}^k$ are projections in $$
For $\{q_i^k\}_{i=1,\cdots,L^k}\in {\rm Proj}
(M_{N_0^k+\cdots+N_k^k+1}(A/I_{k+1}))$ in $(A_k)$
and $\hat P_{k+1}:=\pi_{k+2}(\bar p)\oplus\sum_{i=1}^{L_0^k}\tilde
p_{0i}^k\oplus\cdots \oplus \sum_{i=1}^{L^k_k}\tilde p_{ki}^k
\in{\rm Proj}
(M_{N_0^k+\cdots+N_k^k+1}(A/I_{k+2}))$,
we have 
\[
\sum_{i=1}^{L^k}q_i^k
=\pi_{k+1,k+2}(\hat P_{k+1}).
\]
Therefore, again by Lemma \ref{lem23},
there exist mutually orthogonal
projections $\tilde q_i^k,\;i=1,\cdots,L^k$
in $M_{N_0^k+\cdots+N_k^k+1}(A/I_{k+2})$
such that 
\begin{align}\label{hatpq}
\pi_{k+1,k+2}(\tilde q^{k}_i)=q_i^k,\quad
\sum_{i=1}^{L^k}\tilde q^k_i=\hat P_{k+1},
\end{align}
Now, by (\ref{pif}) in $(A_{k})$,
we have
\begin{align}
\lV
\pi_{k+1,k+2}\lmk
\pi_{k+2}\lmk
\bar p\varphi(f)\bar p\rmk
\oplus
\sum_{i=1}^{L_0^k}f(\xi_{0i}^k)\tilde p^k_{1i}
\oplus\cdots \oplus \sum_{i=1}^{L^k_k}f(\xi_{ki}^k)\tilde p^k_{ki}
-\sum_{i=1}^{L^k}f(\zeta_i^k)\tilde q_i^k
\rmk
\rV<\delta_{k+1},\quad \forall f\in{\cal F}.
\end{align}
By Lemma \ref{lem24}, there exist
$e_i\in{\rm Proj}(\tilde q_i^kM_{N_0^k+\cdots+N_k^k+1}
(I_{k+1}/I_{k+2})\tilde q_i^k)$, $i=1,\cdots,L^k$
and $e:=\sum_{i=1}^{L^k}e_i$ such that
\begin{align}\label{zero}
\lV
\lcm
\pi_{k+2}\lmk
\bar p\varphi(f)\bar p\rmk
\oplus
\sum_{i=1}^{L_0^k}f(\xi_{0i}^k)\tilde p^k_{0i}
\oplus\cdots \oplus \sum_{i=1}^{L^k_k}f(\xi_{ki}^k)\tilde p^k_{ki}
-\sum_{i=1}^{L^k}f(\zeta_i^k)\tilde q_i^k,
e
\rcm
\rV<4\delta_{k+1},
\end{align}
and
\begin{align}\label{ichi}
\lV
(1-e)\lmk
\pi_{k+2}\lmk\bar p\varphi(f)\bar p\rmk
\oplus
\sum_{i=1}^{L_0^k}f(\xi_{0i}^k)\tilde p^k_{0i}
\oplus\cdots \oplus \sum_{i=1}^{L^k_k}f(\xi_{ki}^k)\tilde p^k_{ki}
-\sum_{i=1}^{L^k}f(\zeta_i^k)\tilde q_i^k
\rmk(1-e)
\rV<2\delta_{k+1},
\end{align}
for all $f\in{\cal F}$.
As $e$ commutes with $\tilde q_i^k$s,
(\ref{zero}) means
\begin{align}\label{ce}
\lV
\lcm
\pi_{k+2}\lmk
\bar p\varphi(f)\bar p\rmk
\oplus
\sum_{i=1}^{L_0^k}f(\xi_{0i}^k)\tilde p^k_{0i}
\oplus\cdots \oplus \sum_{i=1}^{L^k_k}f(\xi_{ki}^k)\tilde p^k_{ki},
e
\rcm
\rV<4\delta_{k+1},
\end{align}for all $f\in{\cal F}$.
Using this, (\ref{sct}),
and the fact
\begin{align}\label{eub}
e\le \sum_{i=1}^{L^k}
\tilde q_{i}^k=\hat P_{k+1}=
\pi_{k+2}(\bar p)\oplus\sum_{i=1}^{L_0^k}\tilde
p_{0i}^k\oplus\cdots \oplus \sum_{i=1}^{L^k_k}\tilde p_{ki}^k
\le \pi_{k+2}(\bar p)\oplus 1\oplus\cdots
\oplus 1,
\end{align}
we have
\begin{align}\label{ep}
\lV
\lcm
\pi_{k+2}\lmk
\varphi(f)\rmk
\oplus
\sum_{i=1}^{L_0^k}f(\xi_{0i}^k)\tilde p^k_{0i}
\oplus\cdots \oplus \sum_{i=1}^{L^k_k}f(\xi_{ki}^k)\tilde p^k_{ki},
e
\rcm
\rV
%\nonumber\\
%&\le
%\lV
%\lcm
%\lmk
%\pi_{k+2}\lmk\bar p\rmk\oplus1\oplus
%\cdots\oplus 1
%\rmk
%\lmk
%\pi_{k+2}\lmk
%\varphi(f)\rmk
%\oplus
%\sum_{i=1}^{L_0^k}f(\xi_{0i}^k)\tilde p^k_{0i}
%\oplus\cdots \oplus \sum_{i=1}^{L^k_k}f(\xi_{ki}^k)
%\tilde p^k_{ki}\rmk
%\right.\right.\nonumber\\
%&\left.\left.
%\quad\quad\quad  \times \lmk\pi_{k+2}\lmk\bar p\rmk\oplus1\oplus
%\cdots\oplus 1
%\rmk,
%%\right.\right.\\
%%&\left. \left.
%e
%\rcm
%\rV\nonumber\\
%&+2\lV\lcm
%\lmk
%\pi_{k+2}\lmk\bar p\rmk\oplus 1\oplus
%\cdots\oplus 1
%\rmk,
%\lmk
%\lmk
%\pi_{k+2}\lmk \varphi(f)\rmk \oplus\sum_{i=1}^{L_0^k}
%f(\xi_{0i}^k)\tilde p^k_{0i}
%\oplus\cdots\oplus \sum_{i=1}^{L^k_k}
%f(\xi_{ki}^k)\tilde p_{ki}^k
%\rmk
%\rmk
%\rcm
%\rV\nonumber\\
%&<\lV
%\lcm
%\pi_{k+2}\lmk
%\bar p\varphi(f)\bar p\rmk
%\oplus
%\sum_{i=1}^{L_0^k}f(\xi_{0i}^k)\tilde p^k_{0i}
%\oplus\cdots \oplus \sum_{i=1}^{L^k_k}f(\xi_{ki}^k)\tilde p^k_{ki},
%e
%\rcm
%\rV
%+2\delta\nonumber\\
<
4\delta_{k+1}
+2\delta,
\end{align}
for all $f\in{\cal F}$.
%By the definition, we have
%\begin{align*}
%e\le\sum_{i=1}^{L^k}\tilde q_i^k=\hat P_{k+1}
%\le \pi_{k+2}(\bar p)\oplus 1\oplus\cdots \oplus 1.
%\end{align*}
%
Furthremore, as
$e_i\le \tilde q_i^k,\; i=1,\cdots,L^k$, (\ref{ichi}) means
\begin{align*}
\lV
(1-e)\lmk
\pi_{k+2}\lmk\bar p\varphi(f)\bar p\rmk
\oplus
\sum_{i=1}^{L_0^k}f(\xi_{0i}^k)\tilde p^k_{0i}
\oplus\cdots \oplus \sum_{i=1}^{L^k_k}f(\xi_{ki}^k)\tilde p^k_{ki}
\rmk(1-e)
-\sum_{i=1}^{L^k}f(\zeta_i^k)\lmk
\tilde q_i^k-e_i\rmk
\rV<2\delta_{k+1},
\end{align*}for all $f\in{\cal F}$.\\
Let $\varphi'_{k+2}:C(X_{k+2})\to M_{N_0^k+\cdots+N_k^k+1}(A/I_{k+2})$
be a $*$-homomorphism defined by
\[
\varphi'_{k+2}(g):=\pi_{k+2}\circ\varphi
(\hat g)\oplus\sum_{i=1}^{L_0^k}g(\xi_{0i}^k)
\tilde p_{0i}^k
\oplus\cdots\oplus \sum_{i=1}^{L^k_k}g(\xi_{ki}^k)
\tilde p_{ki}^k,
\]
where $\hat g\in C(X)$ is an extension of $g\in C(X_{k+2})$. 
To see that this is well-defined, let $\hat g_1,\hat g_2\in C(X)$
be two extensions of $g$. Then we have $\hat g_1-\hat g_2\vert_{X_{k+2}}=0$.
By the assumption (\ref{vsp}), 
we have $\pi_{k+2}\circ\varphi(\hat g_1-\hat g_2)=0$.\\
%Furthermore, as $\xi_{ji}^k\in X_{j+1}\subset X_{k+1}\subset X_{k+2}$,
%the $\sum_{i=1}^{L_0^k}g(\xi_{0i}^k)\tilde p_{0i}^k
%\oplus\cdots\oplus \sum_{i=1}^{L^k_l}g(\xi_{ki}^k)
%\tilde p_{ki}^k$ part is well-defined as well.\\
From (\ref{ep}), 
we have
\[
\lV\lcm
e,\varphi_{k+2}'(f\vert_{X_{k+2}})\rcm
\rV
<4\delta_{k+1}+2\delta<
6\delta_{k+1}\le
\delta_D(\frac 13\delta_{k+2},{\cal F}_{k+2}, X_{k+2}),
\]
for all $f\in{\cal F}$.
Therefore, from the condition $D_{{\cal F}_{k+2}}$ of $X_{k+2}$,
%For any unital $C^*$-algebra $\cal B$, unital $*$-homomorphism
%$\varphi :C(X)\to {\cal B}$ and a projection $p\in B$
%with
%\[
%\lV p\varphi(f)-\varphi(f)p\rV<\delta,\quad \forall f\in
%{\cal F},
%\]
%\end{center}
we obtain a positive integer ${N^{k+1}}'=N_{\cal D}(\frac 13\delta_{k+2},{\cal F}_{k+2},X_{k+2})$, points
$\xi_{k+1,i}^{k+1}\in X_{k+2},\; i=1,\cdots, L_{k+1}^{k+1}$,
${\zeta_i^{k+1}}'\in X_{{k+2}},\; i=1,\cdots, {L^{k}}'$, and
two sets of mutually orthogonal projections 
\[
p_{k+1,i}^{k+1}\in {\rm Proj}(M_{{N^{k+1}}'}(eM_{N_0^k+\cdots+N_k^k+1}
(A/I_{k+2})e)),\; i=1,\cdots, L_{k+1}^{k+1}\]
and \[
{q_i^{k+1}}'\in {\rm Proj}(M_{{N^{k+1}}'+1}(eM_{N_0^k+\cdots+N_k^k+1}
(A/I_{k+2})e))\quad i=1,\cdots, L^k.\]
They statisfy
\begin{align}\label{epq}
\sum_{i=1}^{L_{k+1}^{k+1}} p_{k+1,i}^{k+1}=
e\otimes 1_{{N^{k+1}}'},\quad
\sum_{i=1}^{{L^{k}}'} {q_i^{k+1}}'
=e\otimes 1_{({N^{k+1}}'+1)},
\end{align}
and
\begin{align}\label{ni}
\lV
e\varphi_{k+2}'(f\vert_{X_{k+2}})e \oplus\sum_{i=1}^{L_{k+1}^{k+1}}
 f(\xi_{k+1,i}^{k+1})p_{k+1,i}^{k+1}
-\sum_{i=1}^{{L^{k}}'} f({\zeta_i^{k+1}}'){q_i^{k+1}}'
\rV<\frac13\delta_{k+2},
\end{align}
for all $f\in{\cal F}$.
As $e\in {\rm Proj}(M_{N_0^k+\cdots+N_k^k+1}(I_{k+1}/I_{k+2}))$,
we have
$
p_{k+1,i}^{k+1}\in 
%{\rm Proj}(M_{{N^{k+1}}'}(M_{N_0^k+\cdots+N_k^k+1}(I_{k+1}/I_{k+2})))
%=
{\rm Proj}(M_{{N^{k+1}}'(N_0^k+\cdots+N_k^k+1)}(I_{k+1}/I_{k+2}))
$
and 
$
{q_{i}^{k+1}}'\in 
%{\rm Proj}(M_{({N^{k+1}}'+1)}(M_{N_0^k+\cdots+N_k^k+1})(A/I_{k+2}))
%=
{\rm Proj}(M_{{(N^{k+1}}'+1)(N_0^k+\cdots+N_k^k+1)}(A/I_{k+2}))$.
Put
 $N_{j}^{k+1}:=N_j^k,\; L^{k+1}_j:=L^{k}_j$,
and $\xi_{ji}^{k+1}
:=\xi_{ji}^k\in X_{j+1}$,
$p^{k+1}_{ji}:=\tilde p_{ji}^{k}\in
{\rm Proj(M_{N_j^{k+1}}(I_j/I_{k+2}))}$ for
$j=0,\cdots,k$, and
$N_{k+1}^{k+1}:={N^{k+1}}'(N_0^k+\cdots+N_k^k+1)$.
We then define unital $*$-homomorphisms
\begin{align*}
h_j^{k+1}:C(X)\to M_{N_j^{k+1}}(I_j/I_{k+2})\\
h_j^{k+1}(f):=\sum_{i=1}^{L^{k+1}_j}
f(\xi_{ji}^{k+1})p_{ji}^{k+1},
\end{align*}
for $j=0,\cdots,k+1$, and
%Furthermore, we define a $*$-homomorphism
\begin{align*}
&H^{k+1}:C(X)\to M_{N_0^{k+1}+\cdots+N_{k+1}^{k+1}+1}(A/I_{k+2})\\
&H^{k+1}(g):=\sum_{i=1}^{L^k}g(\zeta_i^k)
(\tilde q_i^k-e_i)\oplus 0_{M_{{N^{k+1}}'(N_0^k+\cdots+N_k^k+1)}
(A/I_{k+2})}+
\sum_{i=1}^{{L^k}'} g({\zeta_i^{k+1}}')
{q_i^{k+1}}'=:\sum_{i=1}^{L^{k+1}} 
g({\zeta_i^{k+1}})
{q_i^{k+1}}.
\end{align*}
%This 
%is well-defined 
%because $\tilde q_i^k-e_i\le 1-e$
%and ${q_i^{k+1}}'\le e\otimes 1_{M_{{N^{k+1}}'+1}({\mathbb C})}$.
From (\ref{ichi}), (\ref{ni}), and (\ref{ce}), we obtain
\begin{align*}
\lV
\pi_{k+2}(\bar p\varphi(f)\bar p)
\oplus h_0^{k+1}(f)\oplus\cdots\oplus h^{k+1}_{k+1}(f)
-H^{k+1}(f)
\rV
%\\
%\le
%\lV
%(1-e)\lmk
%\pi_{k+2}(\bar p\varphi(f)\bar p)
%\oplus h_0^{k+1}(f)\oplus\cdots\oplus h^{k+1}_k(f)
%-\sum_{i=1}^{L^k}f(\zeta_i^k)
%(\tilde q_i^k-e_i)
%\rmk
%(1-e)
%\rV\\
%+
%\lV
%e\lmk
%\pi_{k+2}(\bar p\varphi(f)\bar p)
%\oplus h_0^{k+1}(f)\oplus\cdots\oplus h^{k+1}_k(f)
%\rmk
%e\oplus h^{k+1}_{k+1}(f)
%-\sum_{i=1}^{{L^k}'} f({\zeta_i^{k+1}}')
%{q_i^{k+1}}'
%\rV\\
%+2\lV
%\lcm
%\pi_{k+2}(\bar p\varphi(f)\bar p)
%\oplus h_0^{k+1}(f)\oplus\cdots\oplus h^{k+1}_k(f),
%e
%\rcm
%\rV\\
<2\delta_{k+1}+\frac 13\delta_{k+2}
+8\delta_{k+1}<\delta_{k+2},
\end{align*}
for all $f\in{\cal F}$.
%Here, we used (\ref{eub}) to see 
%$e=e(\pi_{k+2}(\bar p)\oplus 1\oplus\cdots\oplus 1)$.
Furthermore, by (\ref{hatpq}) and (\ref{epq}), we have
\begin{align*}
\pi_{k+2}(\bar p)\oplus \sum_{i=1}^{L^{k+1}_0}
p_{0i}^{k+1}\oplus\cdots\oplus
\sum_{i=1}^{L^{k+1}_k}
p_{{k},i}^{k+1}\oplus
\sum_{i=1}^{L^{k+1}_{k+1}}
p_{{k+1},i}^{k+1}
%=\pi_{k+2}(\bar p)\oplus \sum_{i=1}^{L^{k}_0}
%\tilde p_{0i}^{k}\oplus\cdots\oplus
%\sum_{i=1}^{L^{k}_k}
%\tilde p_{{k},i}^{k}\oplus
%\sum_{i=1}^{L^{k+1}_{k+1}}
%p_{{k+1},i}^{k+1}=\hat P_{k+1}\oplus
%\sum_{i=1}^{L^{k+1}_{k+1}}
%p_{{k+1},i}^{k+1}\\
%=\sum_{i=1}^{L^k}\tilde q^k_i\oplus \sum_{i=1}^{L^{k+1}_{k+1}}
%p_{{k+1},i}^{k+1}
%=\sum_{i=1}^{L^k}(\tilde q^k_i-e_i)
%\oplus 0+\sum_{i=1}^{L^k}e_i
%\oplus \sum_{i=1}^{L^{k+1}_{k+1}}
%p_{{k+1},i}^{k+1}\\
%=\sum_{i=1}^{L^k}(\tilde q^k_i-e_i)
%\oplus 0+e\otimes 1_{{N^{k+1}}'+1}\\
%=\sum_{i=1}^{L^k}(\tilde q^k_i-e_i)
%\oplus 0+
%\sum_{i=1}^{{L^k}'}{{q_i^{k+1}}'}
=
\sum_{i=1}^{L^{k+1}} 
{q_i^{k+1}}.
\end{align*}
Hence we obtain $(A_{k+1})$.\\
With the same argument, it can be easily checked that ($A_0$) holds.
$\square$
\section{General Interaction}\label{gen}%(See \cite{Simon} for
%the definition and standard
%techniques on quantum spin systems):
The infinite $\nu$-dimensional quantum spin system with one site algebra 
${M_d}({\mathbb C})$ is given by the UHF $C^*$-algebra
\[
{\mathcal A}_{{\mathbb Z}^{\nu}}:=\overline{
\bigotimes_{{\mathbb Z}^{\nu}}{M_d}({\mathbb C})}^{C^*},
\]
which is the $C^*$- inductive limit of the local algebras
\[
\left\{
{\mathcal A}_{\Lambda}:=\bigotimes_{{\Lambda}}{M_d}
({\mathbb C})\vert
\quad \Lambda\subset{\mathbb Z}^{\nu},\quad
| \Lambda|<\infty
\right\}
.
\]
Here, $| \Lambda|$ denotes the number of points in $\Lambda$.
For each $n$, we denote $\nu$-dimensional
 cube $[-n,n]^{\nu}$ by $\Lambda_n$.
Let $\gamma_j,\; j\in {\mathbb Z}^{\nu}$ 
be the $j$-lattice translation.\\
An interaction is a map $\Phi$ from 
the finite subsets
of ${\mathbb Z}^{\nu}$ into ${\mathcal A}_{{\mathbb Z}^{\nu}}$ such
that $\Phi(X) \in {{\cal A}}_{X}$ 
and $\Phi(X) = \Phi(X)^*$
for any $X \Subset {\mathbb Z}^{\nu}$. 
An interaction $\Phi$ is said to be translation-invariant
if \[
\Phi(X+j)=\gamma_j\lmk
\Phi(X)\rmk,\quad \forall j\in{\mathbb Z}^{\nu},\quad
\forall X\Subset {\mathbb Z}^{\nu}.
\]
%For a subset $\Lambda$ of $\mathbb Z$,
%we denote the surface energy of $\Lambda$ by
%\[
%W_\Phi(\Lambda):=\sum_{X\cap \Lambda\neq\phi,X\cap \Lambda^c\neq\phi}
%\Phi(X).
%\]
A norm of an interaction $\Phi$ is defined by
$\lV\Phi\rV\equiv\sum_{X\ni 0}\lv X\rv^{-1}\lV\Phi(X)\rV$.
\begin{cor}
Let $\Phi_1,\Phi_2,\cdots,\Phi_m$
be a finite set of translation invariant interactions
$\lV \Phi_i\rV<\infty,\; i=1,\cdots,m$,
in the $\nu$-dimensional quantum spin system 
${\mathcal A}_{{\mathbb Z}^{\nu}}$.
For each $i=1,\cdots,m$ and $n\in{\mathbb N}, let$
$H_{i,n}$ be an element in ${\al}_n$ given by
\[
H_{i,n}:=\frac{1}{(2n+1)^\nu}\sum_{I\subset [-n,n]^\nu}
\Phi_i(I)\in{\mathcal A}_{n}.
\]Then there exist 
sequences of selfadjoint elements
$Y_{i,n}\in{\mathcal A}_{n},\;i=1,\cdots, m$ such that
\begin{align*}
\lim_{n\to\infty}\lV
H_{i,n}-Y_{i,n}
\rV=0,\\
\left[Y_{i,n},Y_{j,n}\right]=0,
\quad\forall i,j=1,\cdots,m.
\end{align*}
\end{cor}
{\it Proof}\\
Assume the assertion were false.
Then there exists $\varepsilon>0$
such that 
\begin{align}
\limsup_{n\to\infty}
\lmk\inf\left\{
{\rm max}_{1\le i\le m}
\left\{
\lV
H_{i,n}-X_{i,n}\rV
\right\}
\;:\;X_{i,n}\in {\al}_n,\;
\left[
X_{i,n}, X_{j,n}
\right]=0,\; \forall i,j=1,\cdots,m,
\right\}\rmk
\ge \varepsilon.
\end{align}
For this $\varepsilon$, by the condition
$\lV \Phi_i\rV<\infty$,
we may choose $M\in{\mathbb N}$ large enough 
so that
\begin{align}\label{md}
\sum_{X\ni 0,diam X>\sqrt M}\lv X\rv^{-1}\lV\Phi_i(X)\rV
+\sum_{X\ni 0,X\not \subset \Lambda_{M} }
\lv X\rv^{-1}\lV\Phi_i(X)\rV<\frac{\varepsilon}{2}.
\end{align}
For each $n\in{\mathbb N}$,
define
\begin{align*}
K_{i,n}^M:=\frac{1}{\lv \Lambda_{[\frac nM]}\rv}
\sum_{j\in \Lambda_{[\frac nM]}}\gamma_{Mj}\lmk
H_{i,M}
\rmk\in{\cal A}_n.
\end{align*}
Then for this $K_{i,n}^M\in{\al}_n$, by Theorem \ref{main},
there exist $Y_{in}^M\in {\al}_n,\;i=1,\cdots,m$
such that
\begin{align}\label{yk}
\lim_{n\to\infty}\lV Y_{in}^M-K_{in}^M\rV=0,\quad
\left[
Y_{in}^M,Y_{jn}^M
\right]=0.
\end{align}
%On the other hand, 
%\begin{align*}
%\lV
%K_{in}^M-H_{in}
%\rV\le
%\lV
%\frac{1}{\lv\Lambda_{[\frac nM]}\rv}
%\lmk
%\sum_{j\in \Lambda_{[\frac nM]}}\gamma_{Mj}\lmk
%H_{i,M}\rmk-\lv\Lambda_{[\frac nM]}\rv H_{i,M[\frac nM]}
%\rmk
%\rV\\
%+
%\lv
%1-\frac{\lv\Lambda_{M[\frac nM]}\rv}{\lv\Lambda_{n}\rv}
%\rv
%\lV
%H_{i,M[\frac nM]}
%\rV+
%\lV\frac{\lv\Lambda_{M[\frac nM]}\rv}{\lv\Lambda_{n}\rv}
%H_{i,M[\frac nM]}-H_{in}
%\rV\\
%\le
%\frac{1}{\lv\Lambda_{M[\frac nM]}\rv}
%\sum_{x\in \Lambda_{M[\frac nM]}}
%\sum_{X\ni x,X\not\subset Mj+\Lambda_M}
%\lv X\rv^{-1}\lV \Phi_i(X)\rV
%+
%\lV\Phi_i\rV
%\lv1-\frac{\lv \Lambda_{M[\frac nM]}\rv}{\lv\Lambda_n\rv}\rv\\
%+\frac{1}{\lv\Lambda_n \rv}
%\sum_{x\in \Lambda_n}\sum_{X\ni x X\not\subset\Lambda_M[\frac nM]}
%\lv X\rv^{-1}\lV \Phi_i(X)\rV\\
%\le
%\sum_{X\ni 0,X\not\subset \Lambda_M}
%\lv X\rv^{-1}\lV \Phi_i(X)\rV
%+
%\lV\Phi_i\rV
%\lv1-\frac{\lv\Lambda_{M[\frac nM]}\rv}{\lv\Lambda_n\rv}\rv\\
%+\frac{1}{\lv\Lambda_n \rv}
%\sum_{x\in \Lambda_{M[\frac nM]-\sqrt M}}
%\sum_{X\ni x X\not\subset\Lambda_M[\frac nM]}
%\lv X\rv^{-1}\lV \Phi_i(X)\rV
%+\frac{1}{\lv\Lambda_n\rv}
%\sum_{x\in \Lambda_n\backslash\Lambda_{M[\frac nM]-\sqrt M}}\sum_{X\ni x X\not\subset\Lambda_M[\frac nM]}
%\lv X\rv^{-1}\lV \Phi_i(X)\rV\\
%\le 
%\sum_{X\ni 0,X\not\subset \Lambda_M}
%\lv X\rv^{-1}\lV \Phi_i(X)\rV
%+
%\lV\Phi_i\rV
%\lv1-\frac{\lv \Lambda_{M[\frac nM]}\rv}{\lv\Lambda_n\rv}\rv\\
%+\sum_{X\ni 0diam X>\sqrt M}\lv X\rv^{-1}\lV \Phi_i(X)\rV
%+\frac{\lv\Lambda_n\rv -\lv\Lambda_{M[\frac nM]-\sqrt M}\rv}
%{\lv\Lambda_n \rv}\lV\Phi_i\rV.
%\end{align*}
By the standard argument in spin systems (see \cite{or} for example)
 we have
\begin{align}
\limsup_n\lV
K_{in}^M-H_{in}
\rV\le
\sum_{X\ni 0,diam X>\sqrt M}\lv X\rv^{-1}\lV\Phi_i(X)\rV
+\sum_{X\ni 0,X\not \subset \Lambda_{M} }\lv X\rv^{-1}\lV\Phi_i(X)\rV
<\frac{\varepsilon}{2}.
\end{align}
From this and (\ref{yk}), we have
\[
\limsup_{n\to\infty}
\lV Y_{in}^M-H_{in}\rV<\frac{\varepsilon}{2},\quad
i=1,\cdots,m.
\]
This is a contradiction.$\square$

\end{document}